\documentclass{gtart_a}
\pdfoutput=1
\usepackage{graphicx}

\title{Refined Kirby calculus for integral homology spheres}

\author{Kazuo Habiro}
\givenname{Kazuo}
\surname{Habiro}
\address{Research Institute for Mathematical Sciences\\ 
Kyoto University\\\newline 
Kyoto 606--8502\\Japan}
\email{habiro@kurims.kyoto-u.ac.jp}
\urladdr{}

\volumenumber{10}
\issuenumber{}
\publicationyear{2006}
\papernumber{31}
\lognumber{0690}
\startpage{1285}
\endpage{1317}

\doi{}
\MR{}
\Zbl{}

\arxivreference{math.GT/0509039}



\received{20 December 2005}
\revised{23 June 2006}
\accepted{20 June 2006}
\published{18 September 2006}
\publishedonline{18 September 2006}
\proposed{Colin Rourke}
\seconded{Peter Teichner, Rob Kirby}
\corresponding{}
\editor{}
\version{}

\AtBeginDocument{\let\bar\wbar\let\tilde\wtilde\let\hat\what}
\def\SetFigFont#1#2#3#4#5{\footnotesize}

\DeclareMathOperator{\GL}{GL}

\makeatletter
\def\cnewtheorem#1[#2]#3{\newtheorem{#1}{#3}[section]
\expandafter\let\csname c@#1\endcsname\c@theorem}

\theoremstyle{plain}
\newtheorem{theorem}{Theorem}[section]
\cnewtheorem{lemma}[theorem]{Lemma}
\cnewtheorem{corollary}[theorem]{Corollary}
\cnewtheorem{proposition}[theorem]{Proposition}

\theoremstyle{definition}
\newtheorem{definition}{Definition}

\theoremstyle{remark}
\cnewtheorem{remark}[theorem]{Remark}
\newtheorem*{claim}{Claim}

\makeatother  

\newcommand\xto[1]{{\overset{#1}{\longrightarrow}}}
\newcommand\xyto[2]{{\overset{#1}{\underset{#2}{\longrightarrow}}}}

\newcommand\modD {{\mathsf D}}
\newcommand\modE {{\mathcal{E}}}
\newcommand\modG {{\mathcal{G}}}
\newcommand\modL {{\mathcal L}}
\newcommand\modM {{\mathcal{M}}}
\newcommand\N{{\mathsf N}}
\newcommand\modS {{\mathcal S}}
\newcommand\modZ {{\mathbb Z}}
\newcommand\tL{\tilde L}
\newcommand\modP {{\mathsf P}}
\newcommand\modQ {{\mathsf Q}}
\newcommand\modW {{\mathsf W}}
\newcommand\modp {{\mathsf p}}
\newcommand\modq {{\mathsf q}}
\newcommand\modw {{\mathsf w}}
\newcommand\I{{\mathsf I}}
\newcommand\PP{P}
\newcommand\QQ{Q}
\newcommand\WW{W}
\newcommand\sfP{{\mathsf P}}
\newcommand\sfW{{\mathsf W}}
\newcommand\sfQ{{\mathsf Q}}
\newcommand\commute[2]{{#1}{#2}={#2}{#1}}
\newcommand\ijI{i,j\in \I,\ i\neq j}
\newcommand\tM{\wtilde\modM}
\newcommand\simb{\sim_b}
\newcommand\simeqto{\overset{\simeq}{\longrightarrow}}
\newcommand\ol{\overline}
\newcommand\W[2]{W_{#1,#2}}
\newcommand\Wx[3]{\modW _{#1,#2}^{#3 1}}
\newcommand\hLv{\hat{L}^{\sharp}}

\begin{document}

\begin{asciiabstract}
  A theorem of Kirby states that two framed links in the 3-sphere
  produce orientation-preserving homeomorphic results of surgery if
  they are related by a sequence of stabilization and handle-slide
  moves.
  The purpose of the present paper is twofold: First, we give a
  sufficient condition for a sequence of handle-slides on framed links
  to be able to be replaced with a sequences of algebraically canceling
  pairs of handle-slides.  Then, using the first result, we obtain a
  refinement of Kirby's calculus for integral homology spheres which
  involves only $\pm 1$-framed links with zero linking numbers.
\end{asciiabstract}

\begin{htmlabstract} 
  A theorem of Kirby states that two framed links in the 3&ndash;sphere
  produce orientation-preserving homeomorphic results of surgery if
  they are related by a sequence of stabilization and handle-slide
  moves.
  The purpose of the present paper is twofold: First, we give a
  sufficient condition for a sequence of handle-slides on framed links
  to be able to be replaced with a sequences of algebraically canceling
  pairs of handle-slides.  Then, using the first result, we obtain a
  refinement of Kirby's calculus for integral homology spheres which
  involves only &plusmn;1&ndash;framed links with zero linking numbers.
\end{htmlabstract}

\begin{abstract} 
  A theorem of Kirby states that two framed links in the $3$--sphere
  produce orientation-preserving homeomorphic results of surgery if
  they are related by a sequence of stabilization and handle-slide
  moves.
  The purpose of the present paper is twofold: First, we give a
  sufficient condition for a sequence of handle-slides on framed links
  to be able to be replaced with a sequences of algebraically canceling
  pairs of handle-slides.  Then, using the first result, we obtain a
  refinement of Kirby's calculus for integral homology spheres which
  involves only $\pm 1$--framed links with zero linking numbers.
\end{abstract}

\maketitle
{\leftskip 25pt\rightskip25pt\small\it
This paper is dedicated to Professor Yukio Matsumoto on the
occasion of his sixtieth birthday.\par}

\section{Introduction}

Every closed, connected, oriented $3$--manifold is realized as the
result of surgery along a framed link in the $3$--sphere, 
Lickorish~\cite{Lickorish}, Wallace~\cite{Wallace}.  Kirby's calculus of framed
links~\cite{Kirby} states that two framed links in the $3$--sphere
have orientation-preserving homeomorphic results of surgeries if and
only if these two links are related by a sequence of two kinds of
moves: {\em stabilizations} and {\em handle-slides}.  Thus Kirby's
calculus provides a method to study closed $3$--manifolds through a
study of framed links.  One of the most successful applications of
Kirby's calculus is Reshetikhin and Turaev's definition of quantum
$3$--manifold invariants \cite{Reshetikhin-Turaev}, which is
considered to give a mathematical definition of Witten's Chern--Simons
path integral \cite{Witten}.

Kirby's calculus involves all the framed links in the $3$--sphere,
which represent all the closed, connected, oriented $3$--manifolds.
However, one is sometimes interested in a more special class of
$3$--manifolds, eg, integral homology spheres.  It is natural to
expect that, by restricting our attention to a special class of framed
links which can represent all the $3$--manifolds under consideration,
we would be able to obtain a {\em refinement} of Kirby's calculus of
special framed links involving some special types of moves, and
consequently we would be able to obtain better results than what we
would obtain by using Kirby's calculus directly.

The present paper is intended as the first of a series of papers in
which we study such refinements of Kirby's calculus.  The purpose of
the present paper is twofold: First, we establish a general result
about sequences of handle-slides on framed links, which will be used
as a ``main lemma'' in the series of papers.  Second, we use the
main lemma to obtain a refinement of Kirby's calculus for integral
homology spheres.

\begin{figure}[ht!]
    \begin{center}\input{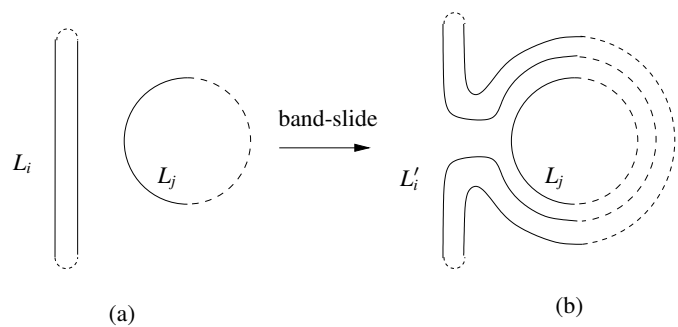}\end{center}
    \caption{(a) Two components $L_i$ and
$L_j$ of a framed link.  (b) The result of a band-slide of $L_i$ over
$L_j$.}
    \label{F10}
  \end{figure}

Let us give a rough description of the main lemma (\fullref{r10}).
Let $M$ be a connected, oriented $3$--manifold, and let $n\ge 0$ be an
integer.  We consider a category $\modS _{M,n}$ whose objects are the
isotopy classes of $n$--component, oriented, ordered, framed links in
$M$, and whose morphisms between two framed links $L$ and $L'$ are
sequences from $L$ to $L'$ of handle-slides, orientation reversals and
permutations.  To each such sequence $S$, we associate in a functorial
way an element $\varphi(S)$ of $\GL(n;\modZ )$, the group of integral
$n\times n$ matrix of determinant~$\pm 1$.  Then the main lemma states that
if the matrix $\varphi(S)$ for $S\co L\to L'$ is the identity matrix
$I_n$, then there is a sequence from $L$ to $L'$ of {\em band-slides}.
A band-slide on a framed link is an algebraically canceling pair of
handle-slides of one component over another, see \fullref{F10}.
Note that if the link is null-homologous in $M$, then a band-slide
preserves the linking matrix.  

It is well known that every integral homology sphere can be expressed
as the result from $S^3$ of surgery along a framed link of diagonal
linking matrix with diagonal entries $\pm 1$.  We call such a framed
link {\em admissible}.  (In the literature, it is also called {\em
algebraically split, unit-framed}.)  Using the main lemma, we can
prove the following refined version of Kirby's calculus for integral
homology spheres.

\begin{theorem}
  \label{t1}
  Two admissible framed links in $S^3$ have orientation-preserving
  homeomorphic results of surgery if and only if they are related by a
  sequence of stabilizations, band-slides and isotopies.
\end{theorem}

Hoste \cite{Hoste} conjectures that if two rationally-framed links in
$S^3$ with zero linking numbers and with framings in $\{1/m\;|\;
m\in \modZ \}$ have orientation-preserving homeomorphic results of
surgery, then they are related by a sequence of Rolfsen's moves
\cite{Rolfsen} through such rationally-framed links.  This conjecture
follows as a corollary to \fullref{t1}, see \fullref{t3}.
We also prove a similar variant of \fullref{t1} for Fenn and
Rourke's theorem \cite{Fenn-Rourke}, see \fullref{t2}.  
\fullref{t1} can also be extended to pairs of integral homology spheres
and knots, see \fullref{r26}, which is a refined version of a
result by Garoufalidis and Kricker~\cite{Garoufalidis-Kricker}.

Now we make some comments on applications of the results in the
present paper.

\begin{remark}
  \label{r19}
  Hoste \cite{Hoste} proves a surgery formula for the Casson invariant
  of integral homology $3$--spheres and shows that if
  \fullref{t3} is true, then his surgery formula provides a
  simple existence proof of the Casson invariant.  This approach to
  the Casson invariant is perhaps the simplest known one if one admits
  \fullref{t3}.
\end{remark}

\begin{remark}
  \label{r28}
  Recall that Ohtsuki's finite type invariants of integral homology
  $3$--spheres \cite{Ohtsuki}, which are generalizations of the Casson
  invariant, are defined in terms of admissible framed links.  Thus it
  is expected that one can use \fullref{t1} in the study of
  Ohtsuki finite type invariants of integral homology spheres.  Though
  this theory of finite type invariants over ${\mathbb Q}$ has been
  understood to a great extent using the Le--Murakami--Ohtsuki invariant
  \cite{Le,LMO}, this is not the case for arbitrary coefficient ring.
  It is expected that, using \fullref{t1}, one can construct a
  universal Ohtsuki finite invariants of integral homology spheres
  over $\mathbb Z$, and perhaps over arbitrary coefficient ring.
\end{remark}

\begin{remark}
  \label{r22}
  In papers in preparation partially joint with T\,T\,Q Le
   \cite{H:in-preparation,Habiro-Le:in-preparation}, we will use
   \fullref{t2} to define, for each simple Lie algebra
   $\mathfrak{g}$, an invariant $J^{\mathfrak g}_M$ of an integral homology
   sphere $M$ which unifies the Witten--Reshetikhin--Turaev invariants
   of $M$ at all roots of unity (for which the invariant is defined),
   which is announced in \cite{H:rims2001}, \cite[Conjecture
   7.29]{Ohtsuki:problem}.  Existence of this invariant implies strong
   integrality properties of the Witten--Reshetikhin--Turaev
   invariants.  \fullref{t2} enables us to prove the
   well-definedness of $J^{\mathfrak g}_M$ without using any previously known
   definitions of the Witten--Reshetikhin--Turaev $3$--manifold
   invariant.  Thus the definition of $J^{\mathfrak g}_M$ provides a new, unified
   definition of the Witten--Reshetikhin--Turaev invariants of integral
   homology spheres.
\end{remark}

We organize the rest of the paper as follows.
In \fullref{sec:defin-stat-main}, we state the main lemma, which
is proved in  \fullref{sec:proof-main-lemma-1}.
In \fullref{sec:proof-theorem-reft1}, we prove \fullref{t1}.
In \fullref{sec:proofs-coroll-reft2}, we prove Hoste's conjectures.
In \fullref{sec:refin-kirby-calc}, we generalize \fullref{t1}
to pairs of integral homology spheres and knots.
In \fullref{sec:remarks-discussions}, we give a short description
of several applications of the main lemma, which we will study in
future papers.

\subsubsection*{Acknowledgements}
  This work started when the the author was a graduate student under
  the direction of Professor Yukio Matsumoto, to whom he would like to
  express his sincere gratitude for continuous encouragement.  He also
  thanks Selman Akbulut, Stavros Garoufalidis, Thang Le, Gregor
  Masbaum, Hitoshi Murakami, Tomotada Ohtsuki and Oleg Viro for
  helpful comments and conversations.
  This research was partially
  supported by the Japan Society for the Promotion of Science,
  Grant-in-Aid for Young Scientists (B), 16740033.

\section{Definitions and the statement of Main Lemma}
\label{sec:defin-stat-main}
In the rest of the paper, all the $3$--manifolds are connected and
oriented.  All homeomorphisms of $3$--manifolds are
orientation-preserving.

In this and the next sections, we fix a connected, oriented
$3$--manifold $M$ and an integer $n\ge 0$.  Let $\modL =\modL _{M,n}$ denote the
set of isotopy classes of $n$--component, oriented, ordered, framed
links in $M$.  We will systematically confuse a framed link and its
isotopy class.  We set $\I=\{1,\ldots ,n\}$.  For $i\in \I$, the $i$th
component of a framed link $L\in \modL $ will be denoted by $L_i$.

\subsection{The category $\modS$ of framed links and elementary moves}
\label{sec:category--framed}

\begin{definition}
  \label{r4}
  Let $\modE =\modE _n$ denote the set of symbols
  \begin{equation*}
    \begin{split}
      \modP _{i,j}&\quad \text{for $\ijI$},\\
      \modQ _i&\quad \text{for $i\in \I$},\\
      \modW _{i,j}^\epsilon &\quad \text{for $\ijI$\ and $\epsilon =\pm 1$}.
    \end{split}
  \end{equation*}
  For $e\in \modE $, an $e$--move on $L\in \modL $ is defined as follows.
  \begin{itemize}
  \item A {\em $\modP _{i,j}$--move} on $L$ exchanges the order of $L_i$
    and $L_j$.
  \item A {\em $\modQ _i$--move} on $L$ reverses the orientation of $L_i$.
  \item A {\em $\modW _{i,j} ^\epsilon $--move} on $L$ is a handle-slide of $L_i$
    over $L_j$.  If $\epsilon =+1$ (resp.\ $\epsilon =-1$), then $L_i$ is added to
    (resp.\ subtracted from) $L_j$, see \fullref{F17}.
  \end{itemize}
  \begin{figure}[ht!]
    \begin{center}\input{\figdir/F17.pstex_t}\end{center}
    \caption{(a) A $\modW _{i,j}^{+1}$--move.  (b) A $\modW _{i,j}^{-1}$--move.}
    \label{F17}
  \end{figure}
  These moves are called {\em elementary moves}.  For $L,L'\in \modL $,
  $e\in \modE $, by $L\xto{e}L'$ we mean that $L'$ is obtained from $L$ by
  an $e$--move.
\end{definition}

If $e=\modP _{i,j}$ or $\modQ _i$, then the result from $L$ of an $e$--move is
unique.  In this case, we denote the result by $e(L)$.  For
$e=\modW _{i,j}^\epsilon $, however, there are in general infinitely many
distinct $L'$ satisfying
$$L\xto{\modW _{i,j}^\epsilon }L'.$$

\begin{definition}
  \label{r1}
  Let $\modS =\modS _{M,n}$ be the free category generated by a
  graph (in the sense of category theory) whose set of vertices are
  $\modL$, and whose edges are elementary moves.  In other words,
  $\modS$ is the category with $\operatorname{Ob}(\modS )=\modL $, and, for
  $L,L'\in \modL $, the set $\modS (L,L')$ of morphisms from $L$ to
  $L'$ consists of the sequences $S=(L^0,e_1,L^1,e_2,L^2,\ldots
  ,e_p,L^p)$ such that $p\ge 0$, $L^0,L^1,\ldots ,L^p\in \modL $,
  $L^0=L$, $L^p=L'$, $e_1,\ldots ,e_p\in \modE $, and for $s=1,\ldots
  ,p$ we have $L^{s-1}\smash{\xto{e_s}}L^s$.  It is convenient to express the
  sequence $S$ as
  \begin{equation*}
    S\co L^0\xto{e_1}L^1\xto{e_2}\cdots\xto{e_p}L^p.
  \end{equation*}
  The identity morphism $1_L\in \modS (L,L)$ of $L\in \modL $ is given by
  $$1_L=(L)\co L\rightarrow L.$$  The composite $S'S$ of
  $S\co L^0\smash{\xto{e_1}}\cdots\smash{\xto{e_p}}L^p$ and
  $S'\co K^0\smash{\xto{e'_1}}\cdots\smash{\xto{e'_p}}K^p$ with $L^p=K^0$
  is given by
  \begin{equation*}
    S'S\co L^0\xto{e_1}\cdots \xto{e_p}L^p=K^0\xto{e'_1}\cdots \xto{e'_p}K^p.
  \end{equation*}
\end{definition}

\subsection{The functor $\varphi\co \modS \rightarrow \GL(n;\modZ )$ and the statement of Main Lemma}
\label{sec:functor-gln-}

For $i,j\in \I$, let $E_{i,j}$ denote the $n\times n$ matrix such that the
$(i,j)$--entry is $1$ and the other entries are $0$.  Let
$I_n=\sum_{i=1}^nE_{i,i}$ be the identity matrix of size $n$.  Define
matrices $\PP_{i,j},\QQ_{i},\WW_{i,j}\in \GL(n;\modZ )$ by
\begin{align*}
  \PP_{i,j}&=I_n-E_{i,i}-E_{j,j}+E_{i,j}+E_{j,i},\\
  \QQ_{i}&=I_n-2E_{i,i},\\
  \WW_{i,j}&=I_n+E_{i,j}
\end{align*}
for $\ijI$.  It is well known that these elements generate
$\GL(n;\modZ )$.  Note that
\begin{equation*}
  \WW_{i,j}^{-1}=I_n-E_{i,j}.
\end{equation*}

We regard the group $\GL(n;\modZ )$ as a category with one object $*$ in
the standard way.
  Define a functor $\varphi\co \modS \to \GL(n;\modZ )$ by $\varphi(L)=*$ for
$L\in \modL $ and
  \begin{equation*}
    \varphi(L\xto{\modP _{i,j}}L')=\PP_{i,j},\quad
    \varphi(L\xto{\modQ _i}L')=\QQ_i,\quad
    \varphi(L\xto{\modW _{i,j}^{\pm 1}}L')=\WW_{i,j}^{\pm 1}.
  \end{equation*}
  For a morphism $S\co L^0\smash{\xto{e_1}}L^1\smash{\xto{e_2}}\cdots
  \smash{\xto{e_p}}L^p$, we have
  \begin{equation*}
    \varphi(S) = \varphi(L^{p-1}\xto{e_p}L^p)\cdots \varphi(L^1\xto{e_2}L^2)\varphi(L^0\xto{e_1}L^1).
  \end{equation*}

Now we state the main lemma.

\begin{theorem}[Main Lemma]
  \label{r10}
  If a morphism $S\co L\rightarrow L'$ in $\modS $ satisfies $\varphi(S)=I_n$, then
  $L$ and $L'$ are related by a sequence of band-slides.
\end{theorem}

\subsection{Linking matrices}
\label{sec:linking-matrices}

  If a framed link $L\in \modL _{M,n}$ is null-homologous in $M$ (ie, each
component of $L$ is null-homologous in $M$), then let $A_L$ denote the
linking matrix of $L$, which is a symmetric matrix with integer
entries of size $n$.  Note that if moreover $S\in \modS (L,L')$,
$L'\in \modL _{M,n}$, then $L'$ also is null-homologous.

For a matrix $T$, let $T^t$ denote the transpose of $T$.

\begin{lemma}
  \label{r9}
  If $L,L'\in \modL _{M,n}$ are null-homologous and $S\in \modS (L,L')$, then we
  have
  \begin{equation*}
    A_{L'}=\varphi(S) A_L \varphi(S)^t.
  \end{equation*}
\end{lemma}

\begin{proof}
  The proof is reduced to the case where $S$ consists of only one
  elementary move, which is well known (see eg Kirby \cite{Kirby:book})
  and can be verified easily.
\end{proof}

\subsection{Explanation using 4--manifolds}
\label{sec:remark-from-4}

The following observation is not necessary in the rest of the paper,
but explains some ideas of the above definitions.

The functor $\varphi\co \modS \rightarrow \GL(n;\modZ )$ has the following natural
topological meaning.  For simplicity, we assume $M=S^3$.  Recall that
for $L\in \modL $, we have a $4$--manifold $X_L$ obtained from the $4$--ball
$B^4$ by attaching $2$--handles $h_1,\ldots ,h_n$ along the components
$L_1,\ldots ,L_n\subset S^3=\partial B^4$ of $L$, see Kirby \cite{Kirby:book}.
The boundary of $X_L$ is the result of surgery $(S^3)_L$.  There is a natural 
basis $u_1,\ldots ,u_n\in H_2(X_L;\modZ )$, where $u_i$ is represented by the 
union of the core of the $2$--handle $h_i$ and the cone of $L_i$ in $B^4$.

Suppose $L\smash{\xto{e}}L'$ with $e\in \modE $.  Then we can define a canonical (up
to isotopy) diffeomorphism $\tilde{e}\co X_L\cong X_{L'}$ as follows.
For $e=\modP _{i,j}$ or $e=\modQ _i$, $\tilde{e}$ is the obvious one.  For
$e=\modW _{i,j}^{\pm 1}$, $\tilde{e}$ is the diffeomorphism given by sliding
$h_i$ along~$h_j$.  Let $u'_1,\ldots ,u'_n$ be the basis of
$H_2(X_{L'};\modZ )$.  Then we have
\begin{equation*}
  \tilde{e}_*=\varphi(e)\co H_2(X_L;\modZ )\rightarrow H_2(X_{L'};\modZ ).
\end{equation*}
Here we regard the matrix $\varphi(e)$ as a $\modZ $--linear map using the
bases of $H_2(X_L;\modZ )$ and $H_2(X_{L'};\modZ )$.
More precisely, we have
\begin{equation*}
  \tilde{e}_*(u_i) = \sum_{j=1}^n \varphi(e)_{i,j} u'_j.
\end{equation*}
For a sequence $S\co L\smash{\xto{e_1}}\cdots \smash{\xto{e_p}}L'$ of elementary moves, the
matrix $\varphi(S)$ corresponds to the isomorphism
$H_2(X_L;\modZ )\rightarrow H_2(X_{L'};\modZ )$ obtained as the composite of the
isomorphisms corresponding to the elementary moves $e_s$, $s=1,\ldots ,p$.

\section[Proof of Theorem~\ref{r10}]{Proof of \fullref{r10}}
\label{sec:proof-main-lemma-1}

\subsection{Bands and annuli for handle-slides}
\label{sec:annuli-handle-1}
In the following, it is sometimes convenient to use bands and annuli
in order to keep track of handle-slides.

\begin{definition}
  \label{r6}
  Suppose that $L\smash{\xto{\modW _{i,j}^\epsilon }}L'$.  \begin{figure}[ht!]
    \begin{center}\input{\figdir/F18.pstex_t}\end{center}
    \caption{}
    \label{F18}
  \end{figure}

  By a {\em band for $L\smash{\xto{\modW _{i,j}^\epsilon}} L'$}, we mean a
  band $b$ in $M$ joining $L_i$ and $L_j$ such that sliding $L_i$ over
  $L_j$ along $b$ (i.e., replacing $L_i$ with a band sum of $L_i$ and a
  parallel copy of $L_j$ along $b$) is a $\modW
  _{i,j}^\epsilon$--move, see \fullref{F18} (a).  In this case, we
  write
$$L\smash{\xyto{\modW _{i,j}^\epsilon }{b}}L'.$$
  By an {\em annulus for $L\smash{\xto{\modW _{i,j}^\epsilon}} L'$}, we mean an annulus
  $a$ in $M$ which looks as depicted in \fullref{F18} (b), such
  that ``handle-slide of $L_i$ over $L_j$ along $a$'' (ie, replacing
  $L_i$ with $L_i'=(L_i\cup \partial a)\setminus \operatorname{int}(L_i\cap \partial a)$) is a
  $\modW _{i,j}^\epsilon $--move.  In this case, we write
$$L\xyto{\modW _{i,j}^\epsilon }{a}L'.$$
\end{definition}

\subsection{Reverse moves and reverse sequences}
\label{sec:reverse-moves-revers}

The {\em reverse} $\wbar{e}$ of $e\in \modE $ is defined by
\begin{gather*}
  \wbar{\modP}_{i,j}=\modP _{i,j},\quad
  \wbar{\modQ}_i=\modQ _i,\quad
  \wwbar{\modW}_{i,j}^\epsilon =\modW _{i,j}^{-\epsilon }.
\end{gather*}

\begin{lemma}
  \label{r7}
  If $e\in \modE $, $L,L'\in \modL $ and $L\xto{e}L'$, then we have
$L'\xto{\wbar{e}}L$.
\end{lemma}

\begin{proof}
  If $e=\modP _{i,j}$ or $e=\modQ _i$, then the assertion is obvious.

  Let $e=\modW _{i,j}^{\epsilon }$.  Choose an annulus $a$ such that
  $L\smash{\xyto{\modW _{i,j}^\epsilon }{a}}L'$.  Then we have
  $L'\smash{\xyto{\modW _{i,j}^{-\epsilon }}{a}}L$.
\end{proof}

For $S\co L^0\smash{\xto{e_1}}L^1\smash{\xto{e_2}}\cdots \smash{\xto{e_p}}L^p$, the {\em reverse}
$\wbar{S}\in \modL (L^p,L^0)$ of $S$ is defined by
\begin{equation*}
  \wbar{S}\co L^p\smash{\xto{\ol{e_p}}}\cdots
\smash{\xto{\ol{e_2}}}L^1\smash{\xto{\ol{e_1}}}L^0.
\end{equation*}
We have
\begin{equation*}
  \varphi(\wbar{S})=\varphi(S)^{-1}.
\end{equation*}

\subsection{Decomposition of $\varphi$}
\label{sec:decomposition--into}

Let $\modM $ denote the free monoid generated by the set $\modE $, which is
regarded as a category with one object $*$.  Define a functor
$$\alpha \co \modS \to \modM $$ by $\alpha (L)=*$ for $L\in \modL $ and
$\alpha (L\smash{\xto{e}}L')=e$ for
$e\in \modE $ with $L\smash{\xto{e}}L'$.    For a morphism
$S\co L^0\smash{\xto{e_1}}L^1\smash{\xto{e_2}}\cdots \smash{\xto{e_p}}L^p$, we have
\begin{equation*}
  \alpha (S) = e_p\cdots e_2e_1.
\end{equation*}
For each element $x=e_p\cdots e_2e_1\in \modM $, the {\em reverse} of $x$ is defined
by
\begin{equation*}
  \wbar{x}= \ol{e_1}\;\ol{e_2}\cdots \ol{e_p}.
\end{equation*}
Clearly, we have $\alpha (\wbar{S})=\ol{\alpha(S)}$ for any morphism $S$ in $\modS $.

Define $\modE ^+$ to be the set of symbols
\begin{equation*}
  \modE ^+=\{\modp _{i,j},\modq _i,\modw _{i,j}\;|\; \ijI\}.
\end{equation*}
Let $\modG $ denote the free group generated by the set $\modE ^+$.  Define a
homomorphism
$$\beta \co \modM \rightarrow \modG $$ by
\begin{gather*}
  \beta (\modP _{i,j})=\modp _{i,j},\quad
  \beta (\modQ _i)=\modq _i,\quad
  \beta (\modW _{i,j}^{\pm 1})=\modw _{i,j}^{\pm 1}.
\end{gather*}
For $x\in \modM $, we have $\beta (\wbar{x})=\beta(x)^{-1}$.

Define a homomorphism $$\gamma \co \modG \to \GL(n;\modZ )$$ by
\begin{equation*}
  \gamma (\modp _{i,j})=\PP_{i,j},\quad
  \gamma (\modq _i)=\QQ_i,\quad
  \gamma (\modw _{i,j})=\WW_{i,j}.
\end{equation*}

Clearly, we have
\begin{equation*}
  \varphi=\gamma \beta \alpha \co \modS  \smash{\xto{\alpha }} \modM
\smash{\xto{\beta }} \modG  \smash{\xto{\gamma }} \GL(n;\modZ ).
\end{equation*}

\subsection{Realization lemma}
\label{sec:realization-lemma}

\begin{lemma}
  \label{r8}
  (1) If $L\in \modL $ and $x\in \modM $, then there is $S\in \modS (L,L')$, $L'\in \modL $,
  such that $\alpha (S)=x$.

  (2) If $L\in \modL $ and $x\in \modM $, then there is $S\in \modS (L',L)$, $L'\in \modL $,
  such that $\alpha (S)=x$.
\end{lemma}

\begin{proof}
  We prove only (1), since (2) can be similarly proved.  Let $l$ be
  the length of $x$.

  If $l=0$, then the result is obvious.

  If $l=1$, then $x\in \modE $.  If $x=\modP _{i,j}$ or $x=\modQ _i$, then set
  $S\co L\smash{\xto{x}} x(L)$.  If $x=\modW _{i,j}^{\epsilon }$, then choose a band $b$
  connecting $L_i$ and $L_j$, and set
$$S\co L\smash{\xyto{\modW _{i,j}^\epsilon }{b}}L',$$
  where $L'$ is the result of $\modW _{i,j}^\epsilon $--move.

  The case $l\ge 2$ reduces to the case $l=1$ by induction.
\end{proof}

\subsection{A preorder on $\modM $}
\label{sec:preorder}
Recall that a {\em preorder} on a set $X$ is a binary relation $\Rightarrow $
such that
\begin{enumerate}
\item $x\Rightarrow x$ for all $x\in X$,
\item $x\Rightarrow y$ and $y\Rightarrow z$ implies $x\Rightarrow z$ for all $x,y,z\in X$.
\end{enumerate}

Define a binary relation $\Rightarrow $ on $\modM $ such that for $x,x'\in \modM $ we have
$x\Rightarrow x'$ if and only if, for any $L,L'\in \modL $ and for any $S\in \modS (L,L')$
with $\alpha (S)=x$, there is $S'\in \modS (L,L')$ satisfying $\alpha (S')=x'$.  (Note
that the definition of $\Rightarrow $ depends on $n$ and $M$.)  It is obvious
that $\Rightarrow $ is a preorder.  By $x\Leftrightarrow y$, we mean ``$x\Rightarrow y$ and $y\Rightarrow x$'',
which is an equivalence relation.

\begin{lemma}
  \label{r15}
For all $x,x',y,y',z\in \modM $, we have the following.
\begin{enumerate}
\item \label{i1} $y\Rightarrow y'$ implies $zyx\Rightarrow zy'x$.
\item \label{i2} $x\Rightarrow y$ implies $\wbar{x}\Rightarrow \wbar{y}$.
\item \label{i4} $yx\Rightarrow z$ implies $x\Rightarrow \wbar{y}z$ and
$y\Rightarrow z\wbar{x}$.
\item \label{i3} $1\Rightarrow \wbar{x} x$.
\end{enumerate}
\end{lemma}

\begin{proof}
  \eqref{i1} Suppose $y\Rightarrow y'$, and $S\co L\rightarrow L'$ with $\alpha (S)=zyx$.  $S$
  can be decomposed into $S=S^3S^2S^1$ with $S^1\co L\rightarrow L^1$,
  $S^2\co L^1\rightarrow L^2$, $S^3\co L^2\rightarrow L'$, $L^1,L^2\in \modL $, such that
  $\alpha (S^1)=x$, $\alpha (S^2)=y$, $\alpha (S^3)=z$.  Since $y\Rightarrow y'$, there is
  $(S^2)'\co L^1\rightarrow L^2$ such that $\alpha ((S^2)')=y'$.  Then we have
  $S^3(S^2)'S^1\co L\rightarrow L'$ and $\alpha (S^3(S^2)'S^1)=zy'x$.  Hence $y\Rightarrow y'$
  implies $zyx\Rightarrow zy'x$.

  \eqref{i2} Suppose $x\Rightarrow y$, and $S\co L\rightarrow L'$ with
  $\alpha (S)=\wbar{x}$.  Then $\wbar{S}\co L'\rightarrow L$ satisfies $\alpha (\wbar{S})=x$.
  Since $x\Rightarrow y$, there is a sequence $S'\co L'\rightarrow L$ such that $\alpha (S')=y$.
  Then we have $\wbar{S}'\co L\rightarrow L'$ and $\alpha (\wbar{S}')=\wbar{y}$.  Hence
  $x\Rightarrow y$ implies $\wbar{x}\Rightarrow \wbar{y}$.

  \eqref{i4} Suppose $yx\Rightarrow z$, and $S\co L\rightarrow L'$ with $\alpha (S)=x$.  By \fullref{r8}, 
  there is a sequence $S'\co L'\rightarrow L''$ with $\alpha (S')=y$.  Since
  $yx\Rightarrow z$, there is a sequence $S''\co L\rightarrow L''$ with $\alpha (S'')=z$.  Then
  we have $\wbar{S}'S''\co L\rightarrow L'$ and $\alpha (\wbar{S}'S'')=\wbar{y}z$.  Hence
  $yx\Rightarrow z$ implies $x\Rightarrow \wbar{y}z$.  Similarly, we can prove that $yx\Rightarrow z$
  implies $y\Rightarrow z\wbar{x}$.

  \eqref{i3} Since $x1=x\Rightarrow x$, it follows from \eqref{i4} that
  $1\Rightarrow \wbar{x}x$.
\end{proof}

\begin{lemma}
  \label{r12}
  We have the following. 
  \begin{enumerate}  \renewcommand{\labelenumi}{(\Alph{enumi})}
  \item \label{A} $\modP _{i,k}\Leftrightarrow \modP _{k,i}$, $\modP _{i,k}^2\Leftrightarrow 1$,
    $\modP _{i,k}\modP _{r,s}\Leftrightarrow \modP _{r,s}\modP _{i,k}$,
    $\modP _{i,k}\modP _{k,r}\Leftrightarrow \modP _{i,r}\modP _{i,k}\Leftrightarrow \modP _{k,r}\modP _{i,r}$,
  \item \label{B} $\modQ _i^2\Leftrightarrow 1$,\qua
    $\modQ _i\modQ _k\Leftrightarrow \modQ _k\modQ _i$,\qua
    $\modQ _j\modP _{i,k}\Leftrightarrow \modP _{i,k}\modQ _j$,\qua
    $\modP _{i,k}\modQ _i\modP _{i,k}\Leftrightarrow \modQ _k$,
  \item \label{C} $\modP _{r,s}\modW _{i,k}^\epsilon \Leftrightarrow \modW _{i,k}^\epsilon \modP _{r,s}$,\qua
    $\modP _{i,k}\modW _{i,k}^\epsilon \Leftrightarrow \modW _{k,i}^\epsilon \modP _{i,k}$,\qua
    $\modP _{i,j} \modW _{i,k}^\epsilon \Leftrightarrow \modW _{j,k}^\epsilon \modP _{i,j} $,

    $\modP _{i,j} \modW _{k,i}^\epsilon \Leftrightarrow \modW _{k,j}^\epsilon \modP _{i,j} $,
  \item \label{D}
    $\modQ _r\modW _{i,k}^\epsilon \Leftrightarrow \modW _{i,k}^\epsilon \modQ _r$,\qua
    $\modQ _k\modW _{i,k}^\epsilon \Leftrightarrow \modW _{i,k}^{-\epsilon }\modQ _k$,\qua
    $\modQ _i\modW _{i,k}^\epsilon \Leftrightarrow \modW _{i,k}^{-\epsilon }\modQ _i$,
  \item \label{E}
    $\modW _{i,k}^\epsilon \modW _{l,m}^\xi \Leftrightarrow \modW _{l,m}^\xi \modW _{i,k}^\epsilon $,\qua
    $\modW _{i,k}^\epsilon \modW _{l,k}^\xi \Leftrightarrow \modW _{l,k}^\xi \modW _{i,k}^\epsilon $,\qua
    $\modW _{i,k}^\epsilon \modW _{i,l}^\xi \Leftrightarrow \modW _{i,l}^\xi \modW _{i,k}^\epsilon $,

    $\modW _{i,k}^\epsilon \modW _{i,k}^\xi \Leftrightarrow \modW _{i,k}^\xi \modW _{i,k}^\epsilon $,
  \item \label{F}
    $\modW _{i,k}^{+1}\Rightarrow \modW _{k,i}^{+1}\modW _{i,k}^{-1}\modP _{i,k}\modQ _i$,
  \item \label{G} $\modW _{i,k}^\xi \modW _{k,l}^\epsilon \Rightarrow \modW _{k,l}^\epsilon \modW _{i,l}^{\epsilon \xi }\modW _{i,k}^\xi $,\qua
    $\modW _{k,l}^\epsilon \modW _{i,k}^\xi \Rightarrow \modW _{i,k}^\xi \modW _{i,l}^{-\epsilon \xi }\modW _{k,l}^\epsilon $,
  \end{enumerate}
  where $i,k$, etc, are distinct elements in $\I$, and
  $\epsilon ,\xi \in \{\pm 1\}$.
\end{lemma}

\fullref{r12} is related to Nielsen's presentation of $\GL(n;\modZ )$.
Taking a look at the statement of \fullref{r14} below before
proceeding may be useful.

\begin{proof}
  \eqref{A}--\eqref{D} are easy and straightforward.

  We will prove \eqref{E}.  They mean that $\modW _{i,k}^\epsilon $ and $\modW _{p,q}^\xi $
  commute up to $\Leftrightarrow $, if $p\neq k$ and $q\neq i$.
  Suppose that
$$L\xyto{\modW _{i,k}^\epsilon }{a}L''\xyto{\modW _{p,q}^\xi }{a'}L',$$
  where $a$ and $a'$ are annuli.  We can move $a'$ by an isotopy of
  $M$ fixing $L''$ as a subset of $M$ so that
  \begin{enumerate}
  \item if $q=k$ (ie, $(p,q)=(l,k),(i,k)$), then $a\cap a'=L''_k$, where
  $L''_k$ denotes the $k$th component of $L''$,
  \item if $q\neq k$ (i.e., $(p,q)=(l,m),(i,l)$), then
  $a\cap a'=\emptyset$.
  \end{enumerate}
  This can be shown as follows.  All the isotopies below fix $L''$ as
  a subset.  First, if $p=i$, then we isotop $a'$ so that $a'\cap L''_i$
  is disjoint from $a$.  Second, if $q=k$, then we isotop $a'$ so that
  in a small neighborhood of $L_k''$, $a$ and $a'$ meet only along
  $L''_k$.  Choose a properly embedded arc $c$ in $a'$ from a point in
  $L''_q$ to a point $L''_p\cap a'$.  Then we isotop $a'$ into a small
  regular neighborhood of $c\cup L''_q$ in $a'$.  Then we can sweep
  $(a\cap a')\setminus L''_q$ out of $a$ by an isotopy.

  We may regard $a'$ as an annulus for a $\modW _{p,q}^\xi $--move
  $$L\smash{\xyto{\modW _{p,q}^\xi }{a'}}L''',$$
  where $L'''\in \modL $. Since $p\neq k$ and
  $q\neq i$, we have
  $$L\smash{\xyto{\modW _{p,q}^\xi }{a'}}
    L'''\smash{\xyto{\modW _{i,j}^\epsilon }{a}}L'.$$
  This shows
  the direction $\Rightarrow $.  The other direction is similar.

  We prove \eqref{F}.  Suppose
$$L \smash{\xyto{\modW _{i,k}^{+1}}{b}} L'.$$
  Let $V$ be a small regular
  neighborhood of $L_i\cup L_k\cup b$, which is a handlebody of genus $2$.
  The inside of $V$ looks as depicted in the upper left corner of
  \fullref{F04}.  \begin{figure}[ht!]
    \begin{center}\input{\figdir/F04.pstex_t}\end{center}
    \caption{}
    \label{F04}
  \end{figure}
  The result of $\modW _{i,k}^{+1}$--move along $b$ is depicted in the
  upper right corner.  There is a sequence
  \begin{equation*}
    L \smash{\xto{\modQ _i}} L^1 \smash{\xto{\modP _{i,k}}} L^2
    \smash{\xyto{\modW _{i,k}^{-1}}{b}}
    L^3 \smash{\xyto{\modW _{k,i}^{+1}}{b'}} L'
  \end{equation*}
  as depicted in \fullref{F04}, which implies \eqref{F}.

  We prove the first formula in \eqref{G} for $\epsilon =\xi =1$.  The second formula
  can be proved similarly.  Suppose
  \begin{equation}
    \label{e2}
    L\smash{\xyto{\modW _{k,l}^{+1}}{a}}L'
     \smash{\xyto{\modW _{i,k}^{+1}}{b}}L'',
  \end{equation}
  where $a$ is an annulus and $b$ is a band.  By moving $b$ with an
  isotopy of $M$ fixing $L'$ as a subset of $M$, we may assume that
  $a$ and $b$ are disjoint.  Let $V$ be a small regular neighborhood
  of $L_i\cup L_k\cup L_l\cup a\cup b$ in $M$, which is a handlebody of genus
  $3$.  The inside of $V$ looks as depicted in the upper left corner
  of \fullref{F05}, where the sequence \eqref{e2} is depicted in
  the top row.  \begin{figure}[ht!]
    \begin{center}\input{\figdir/F05.pstex_t}\end{center}
    \caption{}
    \label{F05}
  \end{figure}

  There is a sequence
  \begin{equation*}
    L
    \smash{\xyto{\modW _{i,k}^{+1}}{b}} L^1
    \smash{\xyto{\modW _{i,l}^{+1}}{b'}} L^2
    \smash{\xyto{\modW _{k,l}^{+1}}{b''}} L'
  \end{equation*}
  as depicted in \fullref{F05}.  Hence we have the first formula
  for $\epsilon =\xi =1$.  The general case of the first formula can be
  obtained by conjugating the formula by~$\modQ _i^\xi \modQ _k^\epsilon $.
\end{proof}

\begin{lemma}
  \label{r21}
  If $S\co L\rightarrow L'$ with $\alpha (S)=\modW _{i,k}^{\mp1}\modW _{i,k}^{\pm 1}$, then $L$
  and $L'$ are related by a band-slide.
\end{lemma}

\begin{proof}
  Consider the case $(p,q)=(i,k)$ in the proof of \eqref{E} of
  \fullref{r12}.  We may assume
$$L\smash{\xyto{\modW _{i,k}^{\pm1}}{a}}
  L''\smash{\xyto{\modW _{i,k}^{\mp1}}{a'}}L',$$
where $a\cap a'=L''_k$.
  Thus $a,a',L_i,L_k$ look as depicted in \fullref{F07} (a).  By
  isotopy, we obtain \fullref{F07} (b).  By a band-slide of $L_i$
  over $L_k$ as indicated in the figure, we obtain a framed link,
  which is isotopic to $L'$.  \begin{figure}[ht!]
    \begin{center}\input{\figdir/F07.pstex_t}\end{center}
    \caption{}
    \label{F07}
  \end{figure}
\end{proof}

\subsection{Realization of relators in $\GL(n;\modZ )$}
\label{sec:presentation-gln-}
The following lemma follows from a result of Nielsen \cite{Nielsen},
see Magnus--Karrass--Solitar \cite[Section 3.5]{MKS}.

\begin{lemma}[Nielsen]
  \label{r14}
  The group $\GL(n;\modZ )$ has a presentation such that the generators are
  the elements of $\modE ^+$ and the relators are as follows.
  \begin{enumerate}  \renewcommand{\labelenumi}{(\alph{enumi})}
  \item \label{a} $\modp _{i,k}=\modp _{k,i}$,\qua $\modp _{i,k}^2=1$,\qua
    $\commute{\modp _{i,k}}{\modp _{r,s}}$,\qua
    $\modp _{i,k}\modp _{k,r}=\modp _{i,r}\modp _{i,k}=\modp _{k,r}\modp _{i,r}$,
  \item \label{b} $\modq _i^2=1$,\qua
    $\commute{\modq _i}{\modq _k}$,\qua
    $\commute{\modq _j}{\modp _{i,k}}$,\qua
    $\modp _{i,k}\modq _i\modp _{i,k}=\modq _k$,
  \item \label{c} $\modp _{r,s}\modw _{i,k}=\modw _{i,k}\modp _{r,s}$,\qua
    $\modp _{i,k}\modw _{i,k}=\modw _{k,i}\modp _{i,k}$,\qua
    $\modp _{i,j} \modw _{i,k}=\modw _{j,k}\modp _{i,j} $,

    $\modp _{i,j} \modw _{k,i}=\modw _{k,j}\modp _{i,j} $,
  \item \label{d} 
    $\modq _r\modw _{i,k}=\modw _{i,k}\modq _r$,\qua
    $\modq _k\modw _{i,k}=\modw _{i,k}^{-1}\modq _k$,\qua
    $\modq _i\modw _{i,k}=\modw _{i,k}^{-1}\modq _i$,
  \item \label{e} $\commute{\modw _{i,k}}{\modw _{l,m}}$,\qua
    $\commute{\modw _{i,k}}{\modw _{l,k}}$,\qua
    $\commute{\modw _{i,k}}{\modw _{i,l}}$,
  \item \label{f} $\modw _{i,k}^{-1}\modw _{k,i}\modw _{i,k}^{-1}=\modq _i\modp _{i,k}$,
  \item \label{g} $\modw _{i,k}\modw _{k,l}\modw _{i,k}^{-1}\modw _{k,l}^{-1}=\modw _{i,l}$.
  \end{enumerate}
  Here $i,k$, etc, denote distinct elements in $\I$.
\end{lemma}

\begin{proof}
  Magnus--Karrass--Solitar \cite[Section 3.5, Theorem N1]{MKS} gives a presentation of the
  automorphism group $\Phi _n$ of a free group of rank $n$, with
  generators $P_{i,j},\sigma _i,U_{i,j},V_{i,j}$ (in the notation of
  \cite{MKS}) for $\ijI$.  This presentation of $\Phi _n$ yields a
  presentation of $\GL(n;\modZ )$ (denoted by $\Lambda _n$ in \cite{MKS}) by
  setting $U_{i,j}=V_{i,j}$.  In our notations, $P_{i,j}$ and $\sigma _i$
  are denoted by $\modp _{i,j}$ and $\modq _i$, respectively, and
  $U_{i,j}=V_{i,j}$ are denoted by~$\modw _{i,j}$.

  Then we easily obtain from the presentation given in \cite{MKS} a
  presentation of $\GL(n;\modZ )$ with a set of generators $\modE ^+$ and a set
  of relations consisting of \eqref{a}--\eqref{f} above and the following.
  \begin{enumerate}  \renewcommand{\labelenumi}{(g\arabic{enumi})}
  \item \label{g1} $\modw _{i,k}\modw _{k,l}^\epsilon \modw _{i,k}^{-1}\modw _{k,l}^{-\epsilon }=\modw _{i,l}^\epsilon
    =\modw _{k,l}^{-\epsilon }\modw _{i,k}\modw _{k,l}^\epsilon \modw _{i,k}^{-1}$,
  \item \label{g2} $\modw _{i,k}^{-1}\modw _{k,l}^\epsilon \modw _{i,k}\modw _{k,l}^{-\epsilon }=\modw _{i,l}^{-\epsilon }
    =\modw _{k,l}^{-\epsilon }\modw _{i,k}^{-1}\modw _{k,l}^\epsilon \modw _{i,k}$.
  \end{enumerate}
  It is easy to see that \eqref{g1} and \eqref{g2} reduces to \eqref{g} modulo the other
  relations.
\end{proof}

For each relation of the form $x=y$ in \fullref{r14}, the element
$x^{-1}y\in \modG $ will be called a {\em relator} of $\GL(n;\modZ )$.

\begin{definition}
  \label{r18}
  Define a map (not a homomorphism) $\lambda \co \modG \rightarrow \modM $ as follows.  For
  $x\in \modG $, let $y_1\cdots y_p$ be the shortest word representing $x$ such
  that for $k=1,\ldots ,p$ we have either $y_k\in \modE ^+$ or $y_k^{-1}\in \modE ^+$.
  Then we set $\lambda (x)=\lambda (y_1)\cdots \lambda (y_p)$, where
  \begin{gather*}
    \lambda (\modp _{i,j}^{\pm 1})=\modP _{i,j},\quad
    \lambda (\modq _i^{\pm 1})=\modQ _i,\quad
    \lambda (\modw _{i,j}^{\pm 1})=\modW _{i,j}^{\pm 1}.
  \end{gather*}
  Clearly, we have $\beta \lambda =\operatorname{id}_{\modG }$.
\end{definition}

\begin{lemma}
  \label{r13}
  If $r$ is a relator of $\GL(n;\modZ )$, then we have $1\Rightarrow \lambda (r)$ and
  $1\Rightarrow \lambda (r^{-1})$.
\end{lemma}

\begin{proof}
  By \fullref{r15} \eqref{i2} and $\lambda (r^{-1})=\ol{\lambda (r)}$, it
  suffices to prove $1\Rightarrow \lambda (r)$.

  The assertion follows easily from \fullref{r12} and
  \fullref{r15}, where the relations \eqref{A}--\eqref{G} in \fullref{r12}
  corresponds to \eqref{a}--\eqref{g} in \fullref{r14}.  For example, we show
  the case of the relation \eqref{g}.  In this case, $r=\modw _{k,l}\modw
  _{i,k}\modw _{k,l}^{-1}\modw _{i,k}^{-1}\modw _{i,l}$.  Hence
  \begin{equation*}
    \lambda (r)=
    \modW _{k,l}^{+1}\modW _{i,k}^{+1}\modW _{k,l}^{-1}\modW _{i,k}^{-1}\modW _{i,l}^{+1}.
  \end{equation*}
  By \fullref{r12} \eqref{G}, we have
  \begin{equation*}
    \modW _{i,k}^{+1}\modW _{k,l}^{+1}\Rightarrow \modW _{k,l}^{+1}\modW _{i,l}^{+1}\modW _{i,k}^{+1}.
  \end{equation*}
  Hence by \fullref{r15} \eqref{i4}, we have
  \begin{equation*}
    1\Rightarrow \modW _{k,l}^{+1}\modW _{i,l}^{+1}\modW _{i,k}^{+1}\modW _{k,l}^{-1}\modW _{i,k}^{-1}.
  \end{equation*}
  Since $\modW _{i,l}^{+1}$ commutes up to $\Leftrightarrow $ with $\modW _{i,k}^{+1}$,
  $\modW _{k,l}^{-1}$, $\modW _{i,k}^{-1}$ by \fullref{r12} \eqref{E}, we obtain
  $1\Rightarrow \lambda (r)$ using \fullref{r15} \eqref{i1}.
\end{proof}

\begin{lemma}
  \label{r17}
  For each $x\in \modM $ with $\gamma \beta (x)=I_n$, there is $x'\in \modM $ such that
  $\beta (x')=\beta (x)$ and $1\Rightarrow x'$.
\end{lemma}

\begin{proof}
  Since $\gamma \beta (x)=I_n$, $\beta (x)$ is contained in the normal subgroup of $\modG $
  generated by the relators of $\GL(n;\modZ )$.  ie, we can express $\beta (x)$
  as
  \begin{equation*}
    \beta (x) = (u_1^{-1} r_1^{\epsilon _1} u_1)\cdots (u_p^{-1} r_p^{\epsilon _p} u_p),
  \end{equation*}
  where $p\ge 0$, and for $s=1,\ldots ,p$, $r_s$ is a relator of $\GL(n;\modZ )$,
  $\epsilon _s=\pm 1$, and $u_s\in \modG $.  Set
  \begin{equation*}
    \begin{split}
      x'
      &=(\lambda (u_1^{-1})\lambda (r_1^{\epsilon _1})\lambda (u_1))\cdots
      (\lambda (u_p^{-1})\lambda (r_p^{\epsilon _p})\lambda ( u_p)) \\
      &=(\ol{\lambda (u_1)}\lambda (r_1^{\epsilon _1})\lambda (u_1))\cdots
      (\ol{\lambda (u_p)}\lambda (r_p^{\epsilon _p})\lambda ( u_p)).
    \end{split}
  \end{equation*}
  Clearly, we have $\beta (x')=\beta (x)$.  For each $s=1,\ldots ,p$, we have
  \begin{equation}
    \label{e1}
  \begin{split}
    1
    &\Rightarrow \ol{\lambda (u_s)}\lambda ( u_s)
    \quad \text{by \fullref{r15} \eqref{i3}}\\
    &=\ol{\lambda (u_s)}1\lambda ( u_s)\\
    &\Rightarrow \ol{\lambda (u_s)}\lambda (r_s^{\epsilon _s})\lambda ( u_s)
    \quad \text{by \fullref{r15} \eqref{i1} and $1\Rightarrow \lambda (r_s^{\epsilon _s})$},
  \end{split}
  \end{equation}
  where $1\Rightarrow \lambda (r_s^{\epsilon _s})$ follows from
  \fullref{r13}.  By \fullref{r15} \eqref{i1} and \eqref{e1} for
  $s=1,\ldots ,p$, it follows that $1\Rightarrow x'$.
\end{proof}

\begin{lemma}
  \label{r2}
  For each $x\in \modM $ with $\gamma \beta (x)=I_n$, there is $x''\in \modM $ 
  such that $\beta (x'')=1$ and $x\Rightarrow x''$.
\end{lemma}

\begin{proof}
  Set $x''=\wbar{x}'x$, where $x'\in \modM $ is as in \fullref{r17}.
  We have $\beta (x'')=\beta (\wbar{x}'x)=\beta (x')^{-1}\beta (x)=1$.
  Moreover, using \fullref{r15}, we have $x=1x\Rightarrow \wbar{x}'x=x''$.
\end{proof}

\subsection{Submonoids $\modM ^0$ and $\tM^0$ of $\modM $}
\label{sec:submonoid-0-}

Let $\modM ^0$ be the submonoid of $\modM $ generated by the elements
$\modW _{i,j}^{-1}\modW _{i,j}^{+1}$ for all $\ijI$.  Let $\tM^0$
denote the submonoid of $\modM $ consisting of the elements $x\in
\modM $ such that $x\Rightarrow y$ for some $y\in \modM ^0$.

\fullref{r12} \eqref{A}, \eqref{B} and \eqref{E} imply that if $e\in \modE $, 
then we have $\wbar{e}e\in \tM^0$.  We will freely use this fact in the proof of 
the following lemma.

\begin{lemma}
  \label{r16}
  $\tM^0$ is invariant under ``conjugation'' in $\modM $.  Ie, if
  $x\in \tM^0$ and $y\in \modM $, then we have $\wbar{y}xy\in \tM_0$.
\end{lemma}

\begin{proof}
  Clearly, we may assume that the length of $y$ is $1$, ie, $y\in \modE $.
  Since $x\in \tM^0$, we have
  \begin{equation*}
    x\Rightarrow  \modW _{i_1,j_1}^{-1}\modW _{i_1,j_1}^{+1}
    \modW _{i_2,j_2}^{-1}\modW _{i_2,j_2}^{+1}\cdots \modW _{i_t,j_t}^{-1}\modW _{i_t,j_t}^{+1}
  \end{equation*}
  where $t\ge 0$ and $i_s,j_s\in \I,i_s\neq j_s$ for $s=1,\ldots ,t$.
  Hence, using \fullref{r15}, we have
  \begin{equation*}
    \begin{split}
      \wbar{y}xy
      &\Rightarrow  \wbar{y}\modW _{i_1,j_1}^{-1}\modW _{i_1,j_1}^{+1}
      \modW _{i_2,j_2}^{-1}\modW _{i_2,j_2}^{+1}\cdots \modW _{i_t,j_t}^{-1}\modW _{i_t,j_t}^{+1}y\\
      &\Rightarrow  \wbar{y}\modW _{i_1,j_1}^{-1}\modW _{i_1,j_1}^{+1}y\wbar{y}
      \modW _{i_2,j_2}^{-1}\modW _{i_2,j_2}^{+1}y\wbar{y}\cdots
      y\wbar{y}\modW _{i_t,j_t}^{-1}\modW _{i_t,j_t}^{+1}y \\
      &= (\wbar{y}\modW _{i_1,j_1}^{-1}\modW _{i_1,j_1}^{+1}y)
      (\wbar{y}
      \modW _{i_2,j_2}^{-1}\modW _{i_2,j_2}^{+1}y)\cdots
      (\wbar{y}\modW _{i_t,j_t}^{-1}\modW _{i_t,j_t}^{+1}y ).
    \end{split}
  \end{equation*}
  It suffices to show that
  $\wbar{y}\modW _{i_s,j_s}^{-1}\modW _{i_s,j_s}^{+1}y\in \tM^0$ for
  $s=1,\ldots ,t$.  Ie, we may assume that $t=1$, and what we have to
  show is the following:
  \begin{claim}
    If $\ijI$, and $y\in \modE $, then we
  have $\wbar{y}\modW _{i,j}^{-1}\modW _{i,j}^{+1}y\in \tM^0$.
  \end{claim}

  First consider the case where we have
  $\modW _{i,j}^{\pm 1}y\Leftrightarrow y\modW _{i',j'}^{\pm \epsilon }$ for some
  $i',j'\in \I,i'\neq j'$ and $\epsilon =\pm 1$.
  We have
  \begin{equation*}
    \wbar{y}\modW _{i,j}^{-1}\modW _{i,j}^{+1}y
    \Leftrightarrow \wbar{y}\modW _{i,j}^{-1}y\modW _{i',j'}^\epsilon
    \Leftrightarrow \wbar{y}y\modW _{i',j'}^{-\epsilon }\modW _{i',j'}^\epsilon .
  \end{equation*}
  If $y=\modP _{p,q}$ or $\modQ _{p}$, then the claim follows from
  $\wbar{y}y\Rightarrow 1$.  If $y=\modW _{p,q}^{+1}$, then the claim
  immediately holds.  If $y=\modW _{p,q}^{-1}$, then \fullref{r12}
  \eqref{E} implies the claim.

  Now consider the other cases.  We have $y=\modW _{p,q}^{\epsilon }$ with either
  $p=j$ or $q=i$ or both, and $\epsilon =\pm 1$.  It suffices to consider the
  following three cases:

  {\bf Case 1}\qua $(p,q)=(j,k)$, $k\neq i,j$.
  We have
  \begin{align*}
      \modW _{j,k}^{-\xi }\modW _{i,j}^{-\epsilon }\modW _{i,j}^\epsilon \modW _{j,k}^\xi
      &\underset{\text{(G)}}{\Rightarrow }\modW _{j,k}^{-\xi }\modW _{i,j}^{-\epsilon }
      \modW _{j,k}^\xi \modW _{i,k}^{\epsilon \xi }\modW _{i,j}^\epsilon
      \underset{\text{(G)}}{\Rightarrow }\modW _{j,k}^{-\xi }
      \modW _{j,k}^\xi \modW _{i,k}^{-\epsilon \xi }\modW _{i,j}^{-\epsilon }
      \modW _{i,k}^{\epsilon \xi }\modW _{i,j}^{\epsilon }\\
      &\underset{\text{(E)}}{\Rightarrow }\modW _{j,k}^{-\xi }
      \modW _{j,k}^\xi \modW _{i,k}^{-\epsilon \xi }
      \modW _{i,k}^{\epsilon \xi }\modW _{i,j}^{-\epsilon }\modW _{i,j}^\epsilon
      \in \tM^0.
  \end{align*}

  {\bf Case 2}\qua $(p,q)=(k,i)$, $k\neq i,j$.
  We have
  \begin{align*}
      \modW _{k,i}^{-\xi }\modW _{i,j}^{-\epsilon }\modW _{i,j}^\epsilon \modW _{k,i}^\xi
      &\underset{\text{(G)}}{\Rightarrow }\modW _{k,i}^{-\xi }\modW _{i,j}^{-\epsilon }\modW _{k,i}^\xi \modW _{k,j}^{-\epsilon \xi }\modW _{i,j}^\epsilon
      \underset{\text{(G)}}{\Rightarrow }\modW _{k,i}^{-\xi }
      \modW _{k,i}^\xi \modW _{k,j}^{\epsilon \xi }\modW _{i,j}^{-\epsilon }
      \modW _{k,j}^{-\epsilon \xi }\modW _{i,j}^{\epsilon }\\
      &\underset{\text{(E)}}{\Rightarrow }\modW _{k,i}^{-\xi }
      \modW _{k,i}^\xi \modW _{k,j}^{\epsilon \xi }
      \modW _{k,j}^{-\epsilon \xi }\modW _{i,j}^{-\epsilon }\modW _{i,j}^\epsilon
      \in \tM^0.
  \end{align*}

  {\bf Case 3}\qua $(p,q)=(j,i)$.
  We have
  \begin{equation*}
      \modW _{j,i}^{-\xi }\modW _{i,j}^{-\epsilon }\modW _{i,j}^\epsilon \modW _{j,i}^\xi
      \underset{\text{(E)}}{\Rightarrow }\modW _{j,i}^{-\xi }\modW _{i,j}^\xi
      \modW _{i,j}^{-\xi }\modW _{j,i}^\xi .
  \end{equation*}
  If $\xi =+1$, then
  \begin{multline*}
    \modW _{j,i}^{-1}\modW _{i,j}^{+1}\modW _{i,j}^{-1}\modW _{j,i}^{+1}
    \underset{\text{(F)}}{\Rightarrow }\modW _{j,i}^{-1}
    (\modW _{j,i}^{+1}\modW _{i,j}^{-1}\modP _{i,j}\modQ _i)
    (\modQ _i\modP _{i,j}\modW _{i,j}^{+1}\modW _{j,i}^{-1})
    \modW _{j,i}^{+1}\\
    \underset{\text{(A),(B)}}{\Rightarrow }\modW _{j,i}^{-1}
    \modW _{j,i}^{+1}\modW _{i,j}^{-1}
    \modW _{i,j}^{+1}\modW _{j,i}^{-1}
    \modW _{j,i}^{+1}
    \in \modM ^0.
  \end{multline*}
  If $\xi =-1$, then
  \begin{multline*}
      \modW _{j,i}^{+1}\modW _{i,j}^{-1}\modW _{i,j}^{+1}\modW _{j,i}^{-1}
      \underset{\text{(F)}}{\Rightarrow }\modW _{j,i}^{+1}
    (\modQ _i\modP _{i,j}\modW _{i,j}^{+1}\modW _{j,i}^{-1})
    (\modW _{j,i}^{+1}\modW _{i,j}^{-1}\modP _{i,j}\modQ _i)
    \modW _{j,i}^{-1}\\
    \underset{\text{(A),(B),(C),(D)}}{\Longrightarrow }\modW _{j,i}^{+1}
    \modW _{j,i}^{-1}\modW _{i,j}^{+1}
    \modW _{i,j}^{-1}\modW _{j,i}^{+1}
    \modW _{j,i}^{-1}
    \in \tM^0.
  \end{multline*}
  This completes the proof of the claim, and hence the lemma.
\end{proof}

\begin{lemma}
  \label{r3}
  If $x\in \modM $ and $\gamma \beta (x)=I_n$, then we have $x\in \tM^0$.
\end{lemma}

\begin{proof}
  By \fullref{r2}, we may assume without loss of generality that
  $\beta (x)=1$.  This implies that there is a sequence
  $x_0=1,x_1,\ldots ,x_p=x\in \modM $ such that for each $s=1,\ldots ,p$, $x_s$ is
  obtained from $s$ by inserting $\wbar{e}e$ with $e\in \modE $, ie, we can
  write $x_{s-1}=y_{s-1}z_{s-1}$ and $x_s=y_{s-1}\wbar{e}ez_{s-1}$.
  Hence, by inserting $\wbar{e}e$, $e\in \modE $, finitely many times into
  $x$, we obtain $x'\in \modM $ with $x\Rightarrow x'$ and
  \begin{equation*}
    x'=(\ol{u_1}u_1)(\ol{u_2}u_2)\cdots (\ol{u_q}u_q),
  \end{equation*}
  where $q\ge 0$, $u_1,\ldots ,u_q\in \modM $.  By \fullref{r16}, it follows that
  $\ol{u_t}u_t\in \tM^0$ for $t=1,\ldots ,q$ (using induction on the
  length of $u_t$).  Hence we have $x'\in \tM^0$.  This and $x\Rightarrow x'$
  imply $x\in \tM^0$.
\end{proof}

\subsection[Proof of Theorem~\ref{r10}]{Proof of \fullref{r10}}
\label{sec:proof-theorem-refr10}
By assumption, we have $\gamma\beta\alpha(S)=\varphi(S)=I_n$ for
$S\co L\to L'$ in $\modS$.  By \fullref{r3}, we have
$\alpha(S)\in\tM^0$.  Hence there is $y\in \modM ^0$ such that
$\alpha(S)\Rightarrow y$.  Hence there is $S'\co L\to L'$ with
$\alpha(S')=y\in\modM^0$.  By \fullref{r21}, it follows that there
is a sequence of band-slides from $L$ to $L'$.

\section[Proof of Theorem~\ref{t1}]{Proof of \fullref{t1}}
\label{sec:proof-theorem-reft1}

\subsection{Definitions and notations}
\label{sec:notations}
In this section, we consider null-homotopic framed links in a fixed
oriented $3$--manifold $M$.  Here a framed link $L$ is said to be
{\em null-homotopic} if every component of $L$ is null-homotopic.
For $p,q\ge 0$, set $I_{p,q}=I_p\oplus(-I_q)$, where $\oplus$ denotes
block sum.  Set
\begin{equation*}
  \modL ^0_{M;\,p,q}=\{L\in \modL _{M,\,p+q}\;|\; A_L=I_{p,q},\ \text{$L$ is
  null-homotopic in $M$}\},
\end{equation*}
where $A_L$ denotes the linking matrix of $L$.  Let
$\modS ^0_{M;\,p,q}$ denote the full subcategory of $\modS _{M,\,p+q}$
such that $\operatorname{Ob}(\modS ^0_{M;\,p,q})=\modL ^0_{M;\,p,q}$.

A component $L_i$ of a framed link $L$ is said to be {\em trivial} if
it bounds a disc which is disjoint from the other components of $L$.
Here the framing of $L_i$ may be arbitrary.

For $0\le p\le p'$ and $0\le q\le q'$, we define a {\em stabilization map}
\begin{equation*}
  \iota _{p',q'}\co \modL ^0_{M;\,p,q}\rightarrow \modL ^0_{M;\,p',q'}
\end{equation*}
as follows.  For $L\in \modL ^0_{M;\,p,q}$, let
\begin{equation*}
  \hat L=\iota _{p',q'}(L)\in \modL ^0_{M;\,p',q'}
\end{equation*}
denote the framed link obtained from $L$ by adjoining $p'-p$ trivial,
$+1$--framed components $O^+_1,\ldots ,O^+_{p'-p}$, and $q'-q$ trivial,
$-1$--framed components $O^-_1,\ldots ,O^-_{q'-q}$, so that
\begin{equation*}
  \hat L = (L_1,\ldots ,L_p,O^+_1,\ldots ,O^+_{p'-p},
  L_{p+1},\ldots ,L_{p+q},O^-_1,\ldots ,O^-_{q'-q}),
\end{equation*}
where we express the ordered link $\hat L$ as a sequence of
components.

By abuse of notation, we extend this $\iota $ notation for elementary move
sequences and matrices.  For $0\le p\le p'$ and $0\le q\le q'$, define a map
\begin{equation*}
  \iota _{p',q'}\co \modS ^0_{M;\,p,q}(L,L')\rightarrow \modS ^0_{M;\,p',q'}(\iota _{p',q'}(L),\iota _{p',q'}(L'))
\end{equation*}
such that for $S\in \modS ^0_{M;\,p,q}(L,L')$, $\iota _{p',q'}(S)$ is defined to be
the obvious sequence of moves from $L$ to $L'$ obtained from $S$ by
adjoining trivial components which are not involved in the sequence of
moves.  The $\iota _{p',q'}$ defines a functor
\begin{equation*}
  \iota _{p',q'}\co \modS ^0_{M;\,p,q}\rightarrow \modS ^0_{M;\,p',q'}.
\end{equation*}
We also define a homomorphism
\begin{equation*}
  \iota _{p',q'}\co \GL(p+q;\modZ )\rightarrow \GL(p'+q';\modZ )
\end{equation*}
as follows.  For a matrix $T=\left(\begin{matrix} T_{++}&T_{+-}\\
T_{-+}&T_{--}
\end{matrix}\right)\in \GL(p+q;\modZ )$ with $\operatorname{size}(T_{++})=p$,
$\operatorname{size}(T_{--})=q$, set
\begin{equation*}
  \iota _{p',q'}(T)=\left(
  \begin{matrix}
    T_{++}& 0 & T_{+-}&0\\
    0& I_{p'-p}& 0 & 0\\
    T_{-+}& 0 & T_{--}&0\\
    0&0&0&I_{q'-q}
  \end{matrix}
  \right).
\end{equation*}
Note that if $S\in \modS ^0_{M;\,p,q}(L,L')$, $L,L'\in \modL ^0_{M;\,p,q}$, then we have
\begin{equation*}
  \varphi(\iota _{p',q'}(S))=\iota _{p',q'}(\varphi(S)).
\end{equation*}
For $L,L'\in \modL ^0_{M;\,p,q}$, by $L\simb L'$ we mean that there is a
sequence from $L$ to $L'$ of isotopies and band-slides.

\subsection[Proof of Theorem~\ref{t1}]{Proof of \fullref{t1}}

\fullref{t1} follows from the case $M=S^3$ of \fullref{r23}
below, which will be proved in the following subsections.

\begin{theorem}
  \label{r23}
  Let $M$ be a connected, oriented $3$--manifold.  Let
  $L,L'\in \modL ^0_{M;\,p,q}$ and suppose that
  $\modS ^0_{M;\,p,q}(L,L')\neq\emptyset$.  Then for some $p'\ge p$, $q'\ge q$,
  we have $\iota _{p',q'}(L)\simb\iota _{p',q'}(L')$.
\end{theorem}

To prove \fullref{t1}, we need only the case $M=S^3$ of \fullref{r23}.  It is for 
later convenience that we state \fullref{r23} in a general form.

\begin{proof}[Proof of \fullref{t1} assuming \fullref{r23}]
  The ``if'' part is obvious.  We prove the ``only if'' part below.

  Suppose that two admissible, unoriented, unordered framed links
  $\tL$ and $\tL'$ in~$S^3$ have homeomorphic results of surgery.  By
  Kirby's theorem, $\tL$ and $\tL'$ are related by a sequence of
  handle-slides after adjoining some trivial $\pm1$--framed
  components.  Thus we may assume without loss of generality that
  $\tL$ and $\tL'$ are related by a sequence of handle-slides.

  We choose orientations and orderings of components to $\tL$ and
  $\tL'$, obtaining an oriented, ordered framed links
  $L,L'\in \modL _{S^3,\,n}$, where $n$ is the number of components of $L$ and
  $L'$.  Here $L$ and $L'$ are chosen so that the linking matrix $A_L$
  of $L$ is $I_{p,q}$ with $p,q\ge 0$, and the linking matrix $A_{L'}$
  of $L'$ is $I_{p',q'}$ with $p',q'\ge 0$.  Since $\tL$ and $\tL'$ are
  related by a sequence of handle-slides, the signatures of $A_L$ and
  $A_{L'}$ are the same.  Hence we have $p=p'$, $q=q'$, and
  $A_L=A_{L'}=I_{p,q}$.

  Since there is a sequence from $\tL$ to $\tL'$ of handle-slides,
  there is a sequence from $L$ to $L'$ of handle-slides, orientation
  change, and permutation of components.  In other words,
  $\modS ^0_{S^3;\,p,q}(L,L')\neq\emptyset$.  By \fullref{r23}, there
  are $p''\ge p$ and $q''\ge q$ such that
  $\iota _{p'',q''}(L)\simb\iota _{p'',q''}(L')$ are related by a sequence of
  band-slides.  Hence we have the assertion.
\end{proof}

\subsection{Realizing a matrix as a sequence between unlinks}
\label{sec:real-sequ-betw}
The rest of this section is devoted to the proof of \fullref{r23}.
Fix a connected, oriented $3$--manifold $M$.  For $p,q\ge 0$, we write
$\modL _{p,q}=\modL _{M;\,p,q}$ and $\modS _{p,q}=\modS _{M;\,p,q}$.

For $p,q\ge 0$, set
\begin{equation*}
  O(p,q;\modZ )=\{T\in \GL(p+q;\modZ )\;|\; T I_{p,q} T^t=I_{p,q}\},
\end{equation*}
which is a subgroup of $\GL(p+q;\modZ )$.

In this subsection, we will prove the following lemma.

\begin{lemma}
  \label{r29}
  Let $p,q\ge 2$ and let $U\in \modL ^0_{B^3;\,p,q}$ be an unlink.  If $T\in O(p,q;\modZ )$
  with $p,q\ge 2$, then there is $S\in \modS ^0_{B^3;\,p,q}(U,U)$ such that
  $\varphi(S)=T$.
\end{lemma}

\fullref{r29} holds for any connected, oriented $3$--manifold $M$
instead of a $3$--ball $B^3$, but we need only the case of $B^3$.

To prove \fullref{r29}, we need a set of generators of $O(p,q;\modZ )$.

\begin{lemma}[Wall {\cite[1.8]{Wall}}]
  \label{r20}
  If $p,q\ge 2$, then $O(p,q;\modZ )$ is generated by the matrices
  \begin{gather}
    \label{e6}
    \begin{split}
      &P_{i,j}\quad \text{for $1\le i<j\le p$ and for $p+1\le i<j\le p+q$},\\
      &Q_i\quad \text{for $1\le i\le p+q$},
    \end{split}
  \end{gather}
  and the matrix $D_{p,q}=\iota _{p,q}(D)\in O(p,q;\modZ )$, where we set
  \begin{equation*}
    D=
    \left(\begin{matrix} 1&1&-1&0\\-1&1&0&1\\-1&0&1&1\\0&1&-1&1
    \end{matrix}\right)
    \in O(2,2;\modZ ).
  \end{equation*}
\end{lemma}

\begin{proof}[Proof of \fullref{r29}]
  It suffices to prove \fullref{r29} when $T$ is each of the
  generators of $O(p,q;\modZ )$ given in \fullref{r20}.  If $T=P_{i,j}$
  or $T=Q_i$, then the assertion follows since $U$ is an unlink.

  Let us consider the case $T=D_{p,q}$.  Without loss of generality we
  may assume that $p=q=2$, since the case $p=q=2$ implies the general
  case via the stabilization map $\iota _{p,q}$.

  The upper left corner of \fullref{F08} depicts $U$.  By
  performing four handle-slides as indicated in the first row in the
  figure, we obtain $L\in \modL _{B^3,4}$.  These four handle-slides are
  realized as $S'\in \modS _{B^3,4}(U,L)$ such that
  \begin{equation*}
    \alpha (S')=\Wx43-\Wx13-\Wx42+\Wx12+.
  \end{equation*}
  Similarly, as depicted in the second row in \fullref{F08}, there
  is $S''\in \modS _{B^3,4}(U,L)$ such that
  \begin{equation*}
    \alpha (S'')=\Wx34-\Wx24-\Wx31+\Wx21+.
  \end{equation*}
  (Note that $S''$ is obtained from $S'$ by a permutation of indices
  $1\leftrightarrow2$, $3\leftrightarrow4$.)  Set
  $S=\wbar{S}''S'\in \modS ^0_{2,2}(U,U)$.  We have
  \begin{equation*}
    \alpha (S)=\ol{\alpha (S'')}\alpha (S')
    = \Wx21- \Wx31- \Wx24+ \Wx34+ \Wx43- \Wx13- \Wx42+ \Wx12+,
  \end{equation*}
  and hence
  $$\varphi(S)
    =\W21^{-1}\W31^{-1}\W24\W34\W43^{-1}\W13^{-1}\W42\W12 =D.\proved$$
  \def\tempI{$\Wx43-\Wx13-\Wx42+\Wx12+$}
  \def\tempII{$\Wx34-\Wx24-\Wx31+\Wx21+$}
  \begin{figure}[ht!]
    \begin{center}\input{\figdir/F08.pstex_t}\end{center}
    \caption{}
    \label{F08}
  \end{figure}
\end{proof}

\subsection{Reordering components}
\label{sec:reord-comp}

Let $L\in \modL ^0_{p,q}$, and $p'\ge 2p$, $q'\ge 2q$.  Let
$L^+=(L_1,\ldots ,L_p)$ (resp.\ $L^-=(L_{p+1},\ldots ,L_{p+q})$) be the sublinks
of $L$ consisting of the $+1$--framed (resp.\ $-1$--framed)
components of $L$.  Set
\begin{equation*}
  \hat L=\iota _{p',q'}(L)
  =(L^+,O^{+,p},O^{+,p'-2p},L^-,O^{-,q},O^{-,q'-2q})\in \modL ^0_{p',q'},
\end{equation*}
where $O^{\pm ,k}$ denotes $k$ trivial components of framings~$\pm1$.
(Here, by abuse of notation, the sequence of sublinks in the right
hand side means a sequence of components.)  We also set
\begin{equation*}
  \hLv
  =(O^{+,p},L^+,O^{+,p'-2p},O^{-,q},L^-,O^{-,q'-2q})\in \modL ^0_{p',q'},
\end{equation*}
which is obtained from $\hat L$ by interchanging $L^+$ and $O^{+,p}$,
and interchanging $L^-$ and $O^{-,q}$.

\begin{lemma}
  \label{r25}
  Let $L\in \modL ^0_{p,q}$.  Then there are integers $p'\ge 2p$ and $q'\ge 2q$
  such that in the above notations we have $\hat L\simb\hLv$.
\end{lemma}

\begin{proof}
  We may assume $p,q\ge 2$ without loss of generality.

  Let $n$ be a sufficiently large integer which will be determined
  later.  Set $p'=2p+n$ and $q'=2q+n$.  Define $\hat
  L^k\in \modL ^0_{p',q'}$ for $k=0,\ldots ,p+q$ inductively by
  \begin{gather*}
    \hat L^k=
    \begin{cases}
      \hat L & \text{if $k=0$},\\
      \modP _{k,p+k}(\hat L^{k-1}) &\text{if $1\le k\le p$},\\
      \modP _{p'+k-p,p'+q+k-p}(\hat L^{k-1}) & \text{if $p+1\le k\le p+q$}.
    \end{cases}
  \end{gather*}
  We have $\hLv=\hat L^{p+q}$.

  It suffices to prove the following:

\begin{claim}
  For each $k=1,\ldots ,p+q$, we have $\hat L^{k-1}\simb \hat L^k$.
\end{claim}

  Note that each permutation move
  involved in the definition of $\hat L^k$ permutes a component in
  $L^+$ or $L^-$ and a trivial component in $O^{+,p}$ or $O^{-,q}$,
  respectively.  We have only to show that such a permutation can be
  realized as a sequence of band-slides.

  For simplicity, we assume $k=1$; the other cases are similar.  Since
  $L_1$ is null-homotopic in $M$, $L_1$ can be unknotted after
  performing finitely many crossing changes of strings of $L_1$.
  Since one can take $n$ to be sufficiently large, we can perform each
  of these crossing changes by a band-slide of $L_1$ over a trivial
  component distinct from $L_{p+1}$, see \fullref{F14}.
  \begin{figure}[ht!]
    \begin{center}\input{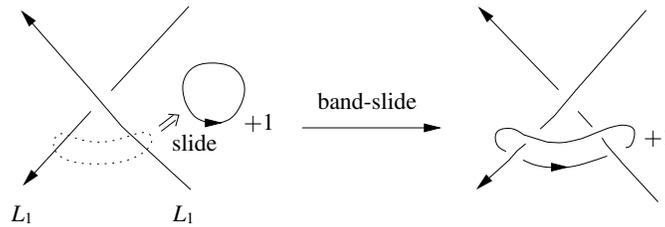}\end{center}
    \caption{Self crossing change realized as a band-slide over
  $+1$--framed trivial component.  The case of the other sign is
  similar.}
    \label{F14}
  \end{figure}  Let $K^1$ denote the framed link obtained from $\hat{L}$
  by applying such band-slides at $c_1,\ldots ,c_r$.  Note that the
  first component $K^1_1$ of $K^1$ is unknotted and of framing~$+1$.
  Let $K^2$ denote the result from $K^1$ by band-sliding all the
  strands linking with $K^1_1$ over $K^1_1$ so that the first
  component $K^2_1$ of $K^2$ is trivial in $K^2$.

  Now we have $K^2=\modP _{1,p+1}(K^2)$ since both the two components
  $K^2_1$ and $K^2_{p+1}$ are trivial in $K^2$.  Hence we have
  \begin{equation*}
    \hat L^0 \simb K^1 \simb K^2
    = \modP _{1,p+1}(K^2)
    \simb \modP _{1,p+1}(K^1)
    \simb \modP _{1,p+1}(\hat L^0)
    =\hat L^1.
  \end{equation*}
  Here the last two $\simb$ can be proved similarly to the first two.
  This completes the proof of the claim, and hence the lemma.  Note
  that it is sufficient to take $n$ as the maximum of the unknotting
  numbers of the components of $L$.
\end{proof}

\subsection{Realizing a matrix as a sequence from one link to itself}
\label{sec:real-gener-case}

\begin{lemma}
  \label{r5}
  Let $L\in \modL ^0_{p,q}$.  There are integers $p'\ge p$, $q'\ge q$ (depending
  on~$L$) such that for each $T\in O(p,q;\modZ )$ there is
  $S\in \modS ^0_{p',q'}(\hat L,\hat L)$, $\hat L=\iota _{p',q'}(L)$, satisfying
  $\varphi(S)=\iota _{p',q'}(T)$.
\end{lemma}

\begin{proof}
  Let $p'\ge p,q'\ge q$, $\hLv\in \modL ^0_{p',q'}$ be as in \fullref{r25}.
  By \fullref{r25}, there is $S'\in \modS ^0_{p',q'}(\hat L,\hLv)$ with
  $\varphi(S')=I_{p'+q'}$.  Note that the sublink
  \begin{equation*}
    L'=\hLv_1\cup \cdots \cup \hLv_p
    \cup \hLv_{p'+1}\cup \cdots \cup \hLv_{p'+q}
  \end{equation*} of $\hLv$ is an unlink
  separated from the other components of $\hLv$ by a sphere.  We can apply
  \fullref{r29} to the sublink $L'$ to obtain a sequence
  $S''\in \modS ^0_{p',q'}(\hLv,\hLv)$ such that
  $\varphi(S'')=\iota _{p',q'}(T)$.
  Set $S=\wbar{S}'S''S'\in \modS ^0_{p',q'}(\hat L,\hat L)$.  Then we have
  \begin{equation*}
    \varphi(S)=\varphi(S')^{-1}\varphi(S'')\varphi(S')
    =\iota _{p',q'}(T).
  \end{equation*}
  This completes the proof.
\end{proof}

\subsection[Proof of Theorem~\ref{r23}]{Proof of \fullref{r23}}
\label{sec:proof-prop-refr23}

Let $S\in \modS ^0_{p,q}(L,L')$.  By \fullref{r5}, there are $p'\ge p$,
$q'\ge q$, $S'\in \modS ^0_{p',q'}(\hat L,\hat L)$, $\hat L=\iota _{p',q'}(L)$ such
that $\varphi(S')=\iota _{p',q'}(\varphi(S))$.  Set
$S''=\iota _{p',q'}(S)\ol{S'}\in \modS ^0_{p',q'}(\hat L,\hat{L}')$,
$\hat{L}'=\iota _{p',q'}(L')$.  Then we have $\varphi(S'')=I_{p'+q'}$.
Hence it follows from \fullref{r10} that $\hat L\simb \hat L'$.  This
completes the proof of \fullref{r23}.

\section{Hoste's conjecture}
\label{sec:proofs-coroll-reft2}

Fenn and Rourke \cite{Fenn-Rourke} prove that Kirby's moves can be
generated by local twisting moves.  Rolfsen \cite{Rolfsen} extends it
to framed links with rational framings.

The purpose of this section is to state and prove ``Fenn--Rourke
version'' and ``Rolfsen version'' of \fullref{t1}, conjectured by
Hoste \cite{Hoste}.

In this section, framed links are unoriented and unordered for
simplicity.

\begin{figure}[ht!]
    \begin{center}\input{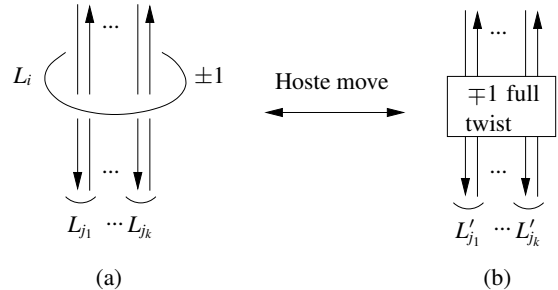}\end{center}
    \caption{A Hoste move.}
    \label{F09}
  \end{figure}

A {\em Hoste move} is defined to be a Fenn--Rourke move between
two admissible framed links, see \fullref{F09}.
  In (a), the component $L_i$ of $L$ is
unknotted and of framing~$\pm1$.  Each of the other components of $L$
links with $L_i$ algebraically $0$ times.  Thus the strands linking
with $L_i$ can be paired as depicted, where $i_1,\ldots ,i_k\neq i$.
(Here $i_1,\ldots ,i_k$ may not be distinct.)  The result $L'$ from
$L$ of a Hoste move on $L_i$ is shown in (b), which is obtained from
$L$ by performing surgery along $L_i$, i.e., by discarding $L_i$ and
giving a $\mp1$ full twist to the bunch of strands linking with $L_i$.

\begin{corollary}[Essentially conjectured by Hoste \cite{Hoste}]
  \label{t2}
   Two admissible framed links in $S^3$ have orientation-preserving
   homeomorphic results of surgery if and only if they are related by
   a sequence of Hoste moves.
\end{corollary}

\begin{figure}[ht!]
    \begin{center}\input{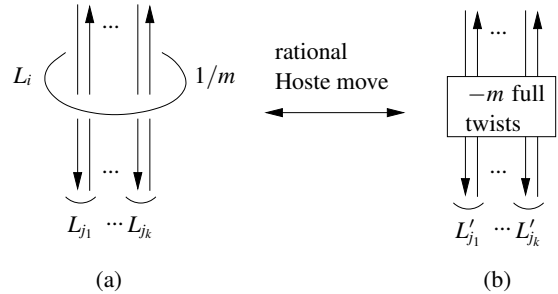}\end{center}
    \caption{A
rational Hoste move.  Here $m$ is any integer.}
    \label{F15}
  \end{figure}

A rationally-framed link in $S^3$ is said to be {\em admissible} if
the linking numbers of any pairs of distinct components are $0$, and
if the framings are in $\{1/m\;|\; m\in \modZ \}$.  Surgery along an admissible
rationally-framed link yields an integral homology sphere.  A {\em
rational Hoste move} is defined to be a Rolfsen move between
two admissible rationally-framed links, see \fullref{F15}.  

\begin{corollary}[Conjectured by Hoste \cite{Hoste}]
  \label{t3}
   Two admissible rationally-framed links in $S^3$ have
   orientation-preserving homeomorphic results of surgery if and only
   if they are related by a sequence of rational Hoste moves.
\end{corollary}

Corollaries~\ref{t2} and~\ref{t3} can be proved by adapting the
proofs by Fenn and Rourke \cite{Fenn-Rourke} and by Rolfsen
\cite{Rolfsen} of the equivalence of their calculi and Kirby's.

\begin{proof}[Proof of \fullref{t2}]
  It is easy to see that a Hoste move can be replaced with a sequence of
  stabilizations and band-slides.  Hence it suffices to prove that a
  band-slide can be replaced with a sequence of Hoste moves.

  Suppose that we are going to perform a band-slide of a component
  $L_i$ of a framed link $L$ over another component $L_j$ of $L$.  By
  finitely many crossing changes for strands in $L_j$ we can unknot
  $L_j$.  This unknotting process can be realized as a sequence of
  finitely many Hoste moves, see \fullref{F12}.  \begin{figure}[ht!]
    \begin{center}\input{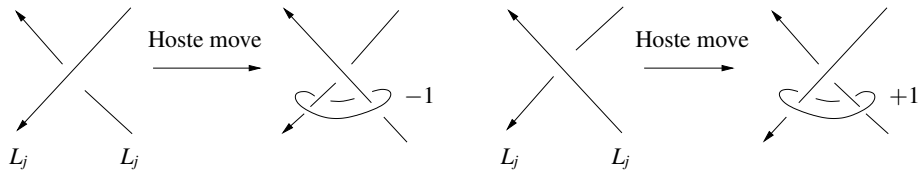}\end{center}
    \caption{A
  realization of crossing change of two strands of $L_j$ as a Hoste
  move.}
    \label{F12}
  \end{figure}  Let $L_j'$ denote the unknotted component obtained from
  $L_j$ by this process, and let $K_1,\ldots ,K_r$ be the newly created
  unknotted components.  A band-slide of $L_i$ over $L_j'$ can then be
  realized by two Hoste moves, see \fullref{F13}.  \begin{figure}[ht!]
    \begin{center}\input{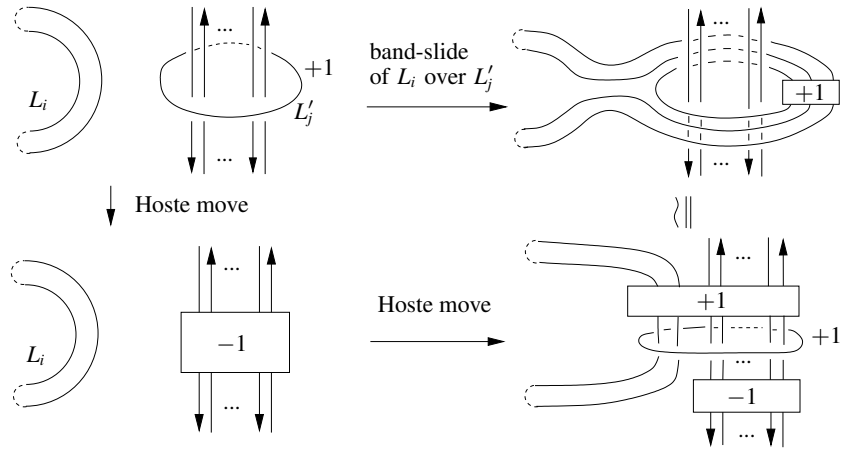}\end{center}
    \caption{The top
  row depicts a band-slide of $L_i$ over a $+1$--framed unknotted
  component~$L_j'$.  This can be replaced with two Hoste moves.  The
  case where $L_j'$ is $-1$--framed is similar.}
    \label{F13}
  \end{figure}  Then we perform
  Hoste moves for the unknotted component $K_1,\ldots ,K_r$.  The result is
  isotopic to the result from $L$ by the band-slide of $L_i$ over
  $L_j$.
\end{proof}

\begin{proof}[Proof of \fullref{t3}]
  It suffices to prove that an admissible rationally-framed link $L$
  is related by a sequence of rational Hoste moves to an admissible
  (integrally) framed link.  The proof is by induction on the number
  of components in $L$ with non-integral framings.  Suppose there is a
  component, say $L_1$, of non-integral framing $1/m$, $m\in \modZ $,
  $m\neq\pm 1$.  We can unknot $L_1$ by some self-crossing changes of
  $L_1$, which can be realized as a sequence of rational Hoste moves
  introducing $\pm 1$--framed components.  Let $L'=L'_1\cup \cdots$ be the
  result of these moves, where $L'_1$ is unknotted and of framing $1/m$.
  Then we perform a rational Hoste move at $L'_1$.  The resulting
  framed link and $L$ have homeomorphic result of surgery.  Moreover, the
  number of components of non-integral framing is reduced by one.
  Hence the assertion follows.
\end{proof}

\section{Knots in integral homology spheres}
\label{sec:refin-kirby-calc}
In this section, framed links are unoriented and unordered for
simplicity.

Let $M$ be a connected, oriented $3$--manifold.
 If a framed link $L$ in $M$ is null-homotopic and if a framed link
$L'$ is related to $L$ by a sequence of Kirby moves, then $L'$ is also
null-homotopic.

A (unoriented, unordered) framed link is said to be {\em
$\pi_1$--admissible} if it is null-homotopic and has diagonal linking
matrix with diagonal entries $\pm 1$.  If $L$ is a null-homotopic
framed link in $L$ with linking matrix of determinant $\pm 1$, then
$L$ is related by a sequence of handle-slides to a $\pi_1$--admissible framed link.

As before, a {\em Hoste move} will mean a Fenn--Rourke move between two
$\pi_1$--admissible framed links.

\begin{proposition}
  \label{r11}
  For two $\pi _1$--admissible framed links $L$ and $L'$ in $M$, the
  following conditions are equivalent.
  \begin{enumerate}
  \item \label{r11.1} $L$ and $L'$ are related by a sequence of Kirby moves (ie,
  stabilizations and handle-slides).
  \item \label{r11.2} $L$ and $L'$ are related by a sequence of stabilizations and
  band-slides.
  \item \label{r11.3} $L$ and $L'$ are related by a sequence of Hoste moves.
  \end{enumerate}
\end{proposition}

\begin{proof}
  Obviously, \eqref{r11.3} implies \eqref{r11.2}, and \eqref{r11.2} implies 
  \eqref{r11.1}.  That \eqref{r11.1} implies \eqref{r11.2} follows easily from 
  \fullref{r23}.  That \eqref{r11.2} implies \eqref{r11.3}
  follows from the proof of \fullref{t2}.  (Note that we need
  the fact that $L$ and $L'$ are null-homotopic, in order to unknot
  some components by Hoste moves in the proofs of \fullref{t2}.)
\end{proof}

Note that any pair $(M,K)$ of an integral homology sphere $M$ and an
oriented knot $K$ in $M$ can be realized as a result from the pair
$(S^3,U)$ of $S^3$ and an unknot $U$ of surgery along a
$\pi _1$--admissible framed link in
$S^3\setminus U$.
Using \fullref{r11}, we have the following refined version of
a theorem by Garoufalidis and Kricker \cite[Theorem
1]{Garoufalidis-Kricker} on surgery presentations of pairs of integral
homology spheres and knots.

\begin{corollary}
  \label{r26}
  Let $L$ and $L'$ be two $\pi _1$--admissible framed links in $S^3\setminus U$.  Then
 the following conditions are equivalent.
  \begin{enumerate}
  \item The results of surgeries, $(S^3,U)_L$ and $(S^3,U)_{L'}$, are
    homeomorphic.
  \item $L$ and $L'$ are related by a sequence of stabilizations and
  band-slides.
  \item $L$ and $L'$ are related by a sequence of Hoste moves.
  \end{enumerate}
\end{corollary}

\begin{proof}
  Garoufalidis and Kricker \cite[Theorem 1]{Garoufalidis-Kricker}
  prove that two null-homotopic framed links in $S^3\setminus U$ with
  linking matrices of determinants $\pm 1$ are related by a sequence
  of Kirby moves if and only if they have homeomorphic results of
  surgeries.  Hence the corollary follows immediately from
  \fullref{r11}.
\end{proof}

\section{Applications}
\label{sec:remarks-discussions}

In this section we describe some applications of Theorems~\ref{r10}
and~\ref{r23}, which we plan to prove in future papers.

\subsection{Splitting the degenerate part}
\label{sec:splitt-link-matr}

A framed link $L$, or the linking matrix $A_L$, is {\em
degenerate-split} if $A_L$ is of the form $O_m\oplus A$, where $m\ge 0$
and $\det A\neq0$.  Note that, for any closed $3$--manifold $M$, there is
a degenerate-split framed link $L$ such that $M\cong S^3_L$.  The
first $m$ components of $L$, which is $0$--framed and has $0$ linking
number with the other components, are called the {\em degenerate
components} or {\em $\modD $--components}.  The other components of $L$ are
called {\em nondegenerate components} or {\em $\N$--components}.  Note
that handle-slide of a component over a $\modD $--component preserves the
linking matrix.

\begin{theorem}
  \label{r24}
  Let $L$ and $L'$ be two (unoriented, unordered) degenerate-split framed
  links in $S^3$.  Then $(S^3_L)\cong(S^3)_{L'}$ if and only if $L$
  and $L'$ are related by a sequence of the following types of moves:
  \begin{itemize}
  \item stabilization, ie, adding or removing a $\pm 1$--framed,
  trivial $\N$--component,
  \item handle-slide of a ($\modD $-- or $\N$--)component over a $\modD $--component,
  \item handle-slide of an $\N$--component over an $\N$--component,
  \item band-slide of a $\modD $--component over an $\N$--component.
  \end{itemize}
\end{theorem}

\begin{remark}
  \label{r31}
  A remarkable application of \fullref{r24} is a refinement of the
  Le--Murakami--Ohtsuki invariant \cite{LMO} of closed, connected,
  oriented $3$--manifolds which is universal for all the
  rational-valued finite type invariants in the sense of Goussarov and
  the author \cite{Goussarov:finite,H}.
\end{remark}

We can also prove the following, which is a generalization of \fullref{t1}.

A framed link $L$ in $S^3$ is {\em split-admissible} if it is
degenerate-split and diagonal with diagonal entries $0,\pm 1$.

\begin{theorem}
  \label{r27}
  Let $L$ and $L'$ be split-admissible framed links in $S^3$.  Then
  $(S^3_L)\cong(S^3)_{L'}$ if and only if $L$ and $L'$ are related by
  a sequence of the following types of moves:
  \begin{itemize}
  \item stabilization,
  \item sliding a ($\modD $-- or $\N$--)component over a $\modD $--component,
  \item band-sliding a ($\modD $-- or $\N$--)component over an $\N$--component.
  \end{itemize}
\end{theorem}

We can also give variants of Theorems~\ref{r24} and~\ref{r27}
involving only local moves, like Fenn and Rourke's theorem or Theorems
\ref{t2} and~\ref{t3}.

It is natural to ask what happens if we drop some of the moves listed
in Theorems~\ref{r24} and~\ref{r27}.  In other words, what kind of
topological structure does the equivalence classes of framed links
correspond to?  For example, we have the following variant of \fullref{r27}.

\begin{theorem}
  \label{r32}
  Let $m\ge 0$.  Let $\modL _m$ denote the set of isotopy classes of
  framed links in $S^3$ with linking matrices of the form $O_m\oplus
  I_{p,q}$, $p,q\ge 0$.  Let $\bar\modL _m$ denote the quotient of $\modL _m$ by
  the equivalence relation generated by stabilization and
  band-sliding.  Let $\modM _m$ denote the set of equivalence classes of
  pairs $(M,f)$ of closed, oriented, spin $3$--manifolds $M$ and an
  isomorphism $f\co \modZ ^m\rightarrow H_1(M;\modZ )$, where two such pairs $(M,f)$ and
  $(M',f')$ are equivalent if there is a spin-structure-preserving
  homeomorphism $\phi \co M\cong M'$ such that $f'=\phi _*f$.  Then there is a
  natural bijection
  \begin{equation*}
    \bar\modL _m \simeqto \modM _m,
  \end{equation*}
  which maps a framed link $L$ to the pair $(S^3_L,f_L)$.  Here the
  spin structure of $S^3_L$ is such that the meridian (with
  $0$--framing in $S^3$) to each $\modD $--component of $L$ represents an
  ``even-framed'' curve in $S^3_L$, and the map $f_L$ maps the $i$th
  basis element of $\modZ ^n$ to the elements represented by the meridian
  (with $0$--framing in~$S^3$) to the $i$th $\modD $--component of $L$.
\end{theorem}

\subsection{Double-slides}
\label{sec:double-slides}
\begin{figure}[ht!]
    \begin{center}\input{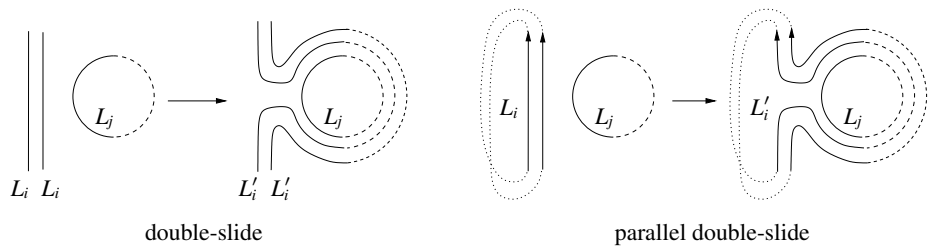}\end{center}
    \caption{A double-slide and a
parallel double-slide.}
    \label{F19}
  \end{figure}

A {\em double-slide} on a framed link $L$ is defined to be
handle-slides of two strands from one component $L_i$ over another
component $L_j$, see \fullref{F19}.   Thus a double-slide is either a band-slide or
a {\em parallel double-slide}, where the two strands are parallel.  It
is easy to see that a band-slide can be realized as a sequence of two
parallel double-slides.

A framed link $L$, or its linking matrix $A_L$, is {\em $2$--diagonal}
if all the non-diagonal entries of $A_L$ are even.  For any symmetric
integer matrix $A$ of size $n$ there is $B\in \GL(n;\modZ )$ such that
$BAB^t$ is $2$--diagonal.  Hence any closed, connected, oriented
$3$--manifold can be obtained from $S^3$ by surgery along a
$2$--diagonal framed link.  A double-slide on a $2$--diagonal framed
link $L$ transforms into another $2$--diagonal framed link, and
preserves the diagonal entries of the linking matrix modulo $4$.

For a $2$--diagonal framed link $L$ with linking matrix $A_L$ with no
diagonal entries congruent to $2$ modulo $4$, define the {\em Brown
number} $b(L)\in \modZ $ of $L$ by
\begin{equation*}
  b(L)=\sigma _4(A_L)-\sigma (A_L).
\end{equation*}
Here $\sigma _4(A_L)=n_1-n_{-1}$, where $n_{\pm 1}$ is the number of diagonal
entries congruent modulo $4$ to $\pm 1$.  $\sigma (A_L)$ denotes the
signature of $A_L$.  If $A_L$ has at least one diagonal entry
$\equiv2\pmod4$, then we formally set $\sigma _4(A_L)=b(L)=\infty $.
$\sigma _4(A_L)\bmod{8}$ is known as the {\em Brown invariant} \cite{Brown}
(see also Matsumoto \cite{Matsumoto}, Kirby--Melvin \cite{KM2}) of the 
$\modZ _4$--valued quadratic form
associated to $A_L$.  Moreover, $b(L)\bmod{8}$ is known to be an invariant
of the $3$--manifold $S^3_L$, called the {\em Brown invariant} of
$S^3_L$, see Kirby--Melvin \cite{KM}.  (Here we formally set $\infty \bmod{8}=\infty $.)  One can
prove that the integer $b(L)$ is invariant under stabilization and
double-slides.  For each $k\in \modZ $, there is a framed link $L^k$ such
that $S^3_{L^k}\cong S^3$ and $b(L^k)=8k$.

A component of a $2$--diagonal framed link is {\em even}
(resp.\ {\em odd}) if its framing is even (resp.\ odd).

We have the following $\modZ _2$--version of \fullref{r27}.

\begin{theorem}
  \label{r34}
  Let $L$ and $L'$ be two $2$--diagonal framed link of the same Brown
  number $n\in \modZ \cup \{\infty \}$.  Then $S^3_L\cong S^3_{L'}$ if and only if
  $L$ and $L'$ are related by a sequence of the following types of
  moves:
  \begin{itemize}
  \item stabilizations,
  \item double-slides,
  \item handle-slides of (even or odd) components over even components.
  \end{itemize}
\end{theorem}

We can modify \fullref{r34} as follows, which may be regarded as
the $\modZ _2$--version of \fullref{r32}.

\begin{theorem}
  \label{r35}
  Let $n\in \modZ \cup \{\infty \}$.  There is a natural bijection between the set
  of $2$--diagonal, oriented, ordered framed links of Brown number $n$
  modulo stabilization and double-slides, and the set of the closed
  $3$--manifolds $M$ of Brown invariant $n\bmod{8}$, equipped with spin
  structure and parameterization of $H_1(M;\modZ _2)$.  The bijection is
  defined similarly as in \fullref{r32}.
\end{theorem}

One can also derive ``local move versions'' of Theorems~\ref{r34} and
\ref{r35}.

\bibliographystyle{gtart}
\bibliography{link}

\end{document}